\documentclass{usiinftr}

\usepackage{tikz}
\usepackage{listings}
\usepackage{amsmath}

\usepackage{booktabs}

\usepackage{subcaption}

\usepackage{siunitx}

\usepackage{color}
\usepackage{tikz}
\usepackage{pgfplots}
\pgfplotsset{compat=1.10}
\usetikzlibrary{patterns}
\usepackage{pgfplotstable}
\usepackage{csvsimple}
\definecolor{excelblue}{RGB}{94,156,211}
\definecolor{excelorange}{RGB}{235,125,60}
\definecolor{excelgray}{RGB}{165,165,165}
\definecolor{excelgreen}{RGB}{114,172,77}

\usepackage{makecell}

\usepackage{xcolor,colortbl}
\definecolor{Gray}{gray}{0.9}
\newcolumntype{G}{>{\columncolor{Gray}}r}

\definecolor{tablecolor1}{HTML}{F0F0F0}
\definecolor{tablecolor2}{HTML}{FFFFFF}

\definecolor{mycolor0}{HTML}{2B83BA}
\definecolor{mycolor1}{HTML}{D7301F}
\definecolor{mycolor2}{HTML}{FC8D59}
\definecolor{mycolor3}{HTML}{abdda4}
\definecolor{mycolor4}{HTML}{e9a3c9}
\definecolor{mycolor5}{HTML}{FEF0D9}

\usepackage{algorithm}
\usepackage{algorithmicx}
\usepackage{algpseudocode}

\newcommand{\Set}[1]{\mathbb{#1}}
\newcommand{\R}[1]{\Set{R}^{#1}}
\newcommand{\ComplexSet}[1]{\Set{C}^{#1}}

\def\bfc{{\discr{c}}}

\def\bfg{{\mathbf g}}
\def\bfh{{\mathbf h}}

\def\bfp{{\mathbf p}}
\def\bfq{{\mathbf q}}

\def\bfu{{\discr{u}}}
\def\bfv{{\mathbf v}}

\def\bfx{{\discr{x}}}

\def\bfI{{\discr{I}}}

\def\F{{\mathbf {F}}}

\def\LineFlow{{\discr{h}}}

\def\PowerFlowSCOPF{{\discr{g}^c}}
\def\LineFlowSCOPF{{\discr{h}^c}}
\def\PowerFlowOPF{{\discr{g}}}
\def\LineFlowOPF{{\discr{h}}}

\newcommand{\pes}[2]{{\mathbf p}^{#1}_{#2}}

\newcommand{\Complex}[1]{\underline{#1}}
\def\Ybc{\Complex{\mathbf{Y}}^{\text{B}}}
\def\vc{\Complex{\mathbf{v}}}
\def\ic{\Complex{\mathbf{I}}}
\def\sc{\Complex{\mathbf{S}}}
\def\sdc{\Complex{\mathbf{S}}^d}
\def\sgc{\Complex{\mathbf{S}}^g}
\def\sfc{\Complex{\mathbf{S}}^f}
\def\stc{\Complex{\mathbf{S}}^t}

\def\bftheta{{\mathbf \theta}}

\def\bflambda{{\boldsymbol \lambda}}

\def\0{{\mathbf 0}}

\def\Nb{{N_\text{B}}}
\def\Ng{{N_\text{G}}}

\def\Nl{{N_\text{L}}}

\newcommand \discr[1]{{\boldsymbol{#1}}}   

\newcommand \Min[1]{{#1}^{\text{min}}}
\newcommand \Max[1]{{#1}^{\text{max}}}

\def\superstar{^{\raise 0.5pt\hbox{$\nthinsp *$}}}

\def\nthinsp{\mskip -2   mu}
\makeatletter
\newcommand*{\transpose}{%
	{\mathpalette\@transpose{}}%
}
\newcommand*{\@transpose}[2]{%
	\raisebox{\depth}{$\m@th#1\intercal$}%
}
\makeatother

\def\minim{\mathop{\hbox{\rm minimize}}}
\def\subject{\text{\rm subject to}}
\def\minimize#1{{\displaystyle\minim_{#1}}}

\begin{document}

\title{Reduced-Space Interior Point \mbox{Methods} in Power Grid Problems} 

\author{Juraj Kardo\v s}{1}
\author{Drosos Kourounis}{1}
\author{Olaf Schenk}{1}

\affiliation{1}{Institute of Computational Science, Universit\`a della Svizzera italiana, Switzerland}

%
%
\TRnumber{2020-1}

%
%

\maketitle

\begin{abstract}
Due to critical environmental issues, the power systems have to accommodate a significant level of penetration of renewable generation which requires smart approaches to the power grid control. 
Associated optimal control problems are large-scale nonlinear optimization problems with up to hundreds of millions of variables and constraints. The interior point methods become computationally intractable, mainly due to the solution of large linear systems. 

This document addresses the computational bottlenecks of the interior point method during the solution of the security constrained optimal power flow problems by applying reduced space quasi-Newton IPM, which could utilize high-performance computers due to the inherent parallelism in the adjoint method. 
Reduced space IPM approach and the adjoint method is a novel approach when it comes to solving the (SC)OPF problems. These were previously used in the PDE-constrained optimization. The presented methodology is suitable for high-performance architectures due to inherent parallelism in the adjoint method during the gradient evaluation, since the individual contingency scenarios are modeled by independent set of the constraints.
Preliminary evaluation of the performance and convergence is performed to study the reduced space approach.    
\end{abstract}

\section{Optimization Problems in the Power Grid}

An electrical grid, or power grid, is an interconnected network for delivering electricity from producers to consumers. It consists of generating stations that produce electrical power, high voltage transmission lines that carry power from distant sources to demand centers and distribution lines that connect individual customers.
More formally, consider a power grid with $\Nb$ buses, $\Ng$ generators, and $\Nl$ transmission lines. The bus voltage vector $\vc \in \ComplexSet{\Nb}$ is defined in polar notation as $\vc = \bfv e ^{j\bftheta}$, where  $\bfv, \bftheta \in \R{\Nb}$ specify the magnitude and phase of the complex voltage. The complex voltages $\vc$ determine the entire power flow (PF) in the grid that can be computed using the Kirchhoff equations and the grid configuration, such as the transmission line parameters, transformer tap ratios, and shunt elements \cite{matpowerManual, meier2006book}. The current injections $\ic \in \ComplexSet{\Nb}$ into the buses are defined as $\ic = \Ybc \vc$, where $\Ybc \in \ComplexSet{\Nb \times \Nb}$ is the bus admittance matrix. The complex power at each bus of the network $\sc = \vc \ic^*$, $\sc \in \ComplexSet{\Nb}$ is  to be balanced by the net power injections from the generators $\sgc \in \ComplexSet{\Ng}$ and demand centers' power consumption $\sdc \in \ComplexSet{\Nb}$. Thus, the alternating current (AC) nodal PF balance equations, also known as the mismatch equations, are expressed as a function of the complex bus voltages and generator injections as $ \sc + \sdc - C_g\sgc= \0$, where $C_g \in \R{\Nb \times \Ng}$ is the generator connectivity matrix, which specifies location of the individual generators in the power grid.

Generator power injections $\sgc = \bfp + j \bfq$ are expressed in terms of active and reactive power components $\bfp, \bfq\in \R{\Ng}$, respectively. Each bus has an associated complex power demand $\sdc$, which is assumed to be known at all of the buses and is modeled by a static ZIP model~\cite{matpowerManual}. If there are no loads connected to bus $i$, then $\{\sdc\}_i = \0$.

The real-world transmission lines are limited by the instantaneous amount of power that can flow through the lines due to the thermal limits. The apparent PFs of the transmission lines, $\sfc \in \ComplexSet{\Nl}$ and $\stc \in \ComplexSet{\Nl}$, are therefore limited by the power injections at both ``from" and ``to" ends of the line, which cannot exceed prescribed limits $\Max{\F}_\text{L}$. Squared values of the apparent PF are usually used in practice.

\section{Optimal Power Flow}

Since the formulation of optimal power flow (OPF) by Carpentier~\cite{Carpentier:1962} as a continuous nonlinear programming (NLP) problem, OPF has become one of the most important and widely studied constrained nonlinear optimal control problems. It is concerned with the optimization of the operation of an electric power grid subject to physical constraints imposed by electrical laws and engineering limits. The objective is to identify the operating configuration that best meets a particular set of evaluation criteria. These criteria may include the cost of generation, transmission line losses, and various requirements concerning the system's security, or resilience with respect to disturbances. The standard OPF problem is formulated as minimization of the active power generation cost, subject to a set of equality and inequality constraints:

\begin{subequations}
	\label{OPF}
	\begin{align}
	\label{OPF:objective}
	\minimize{\bftheta,\bfv,\bfp,\bfq}  & \sum_{l=1}^{\Ng} a_l(\pes{l}{})^2 + b_l\pes{l}{} +  c_l \\ 
	\subject{\,\;}  \label{OPF:PowerFlow}
	& \PowerFlowOPF(\vc, \sgc) = \0, \\  
	\label{OPF:LineFlowLimits}
	& \LineFlowOPF(\vc)  \leq \0, \\
	\label{OPF:VoltageLimits}
	& \Min{\bfv}      \leq \bfv      \leq \Max{\bfv},\quad \bftheta^{ref} = 0,  \\
	\label{OPF:PowerLimits}
	& \Min{\bfp} \leq \bfp  \leq \Max{\bfp},  \quad 
	\Min{\bfq} \leq \bfq  \leq \Max{\bfq}.
	\end{align}
\end{subequations}
The objective function is a quadratic functional of the active power generation, with the cost coefficients $a,b,c$ defining the cost function of each generator. The AC nodal PF balance equations for each bus in the grid represent the nonlinear equality constraints \eqref{OPF:PowerFlow}, $\PowerFlowOPF \in \mathbb{R}^{2\Nb}$, where the exponential terms in the polar voltage representation are replaced using the Euler identity, and subsequently, the complex quantities are split into their real and imaginary components, $\PowerFlowOPF = (\Re(\bfg),\ \Im(\bfg))$. Line power flows limits $\LineFlowOPF \in \mathbb{R}^{2\Nl}$ are imposed by a set of inequality constraints \eqref{OPF:LineFlowLimits}, where 
\begin{equation}
\LineFlowOPF = \left(\sfc(\sfc)^* - (\Max{\F}_\text{L})^2,\quad  \stc(\stc)^* - (\Max{\F}_\text{L})^2 \right).    
\end{equation}

Over the last five decades almost every mathematical programming approach that can be applied to OPF has been attempted. Sequential linear programming (LP) has been presented by Stott and Hobson \cite{Stott:1978I} and Alsac et al.~\cite{Alsac:1990}. Despite its performance it introduces inaccuracies and might lead to infeasible solutions. Sequential quadratic programming (SQP) was applied by Burchett~\cite{Burchett:1982} delivering accurate and feasible solutions but the computational times are relatively high, rendering the whole solution approach inappropriate for large-scale power systems.  
IPMs  have became a successful tool for solving constrained optimization problems. The origin of the IPMs popularity reaches back to the 1984 when Karmarkar \cite{Karmarkar1984} announced a polynomial time linear program that was considerably faster than the most popular simplex method to date. Furthermore, IPMs can be applied also to quadratic and other nonlinear programs, unlike the simplex method which can be applied only to linear programming. 
The main advantages of the IPM lie in the ease of handling inequality constraints by logarithmic barrier functions, speed of convergence, and  a fact that strictly feasible initial point is not required. Linear and nonlinear IPM were applied to the OPF problems since the 90s by Vargas et al.~\cite{Vargas:1993,Granville:1994,Torres:1998}. Several extensions of these algorithmic approaches are reported in recent papers; see, e.g.,~\cite{Frank2012a, Duan:2015, Lu:2018}. Although the LP-based approaches can already be applied close to real time for various practical OPF problems at least in medium size systems, further improvement in tractability need to be envisaged to get closer to the real-time AC OPF dream in large-scale systems, as pointed out by \cite{CAPITANESCU:2016}. Another advantage of IPMs is that they are applicable to large problems and can easily exploit problem structure. The underlying linear algebra components of structured problems allow for large-scale parallel optimization problems on high-performance computers.

Modern trends in power grid operations and modelling render approximation based optimization techniques less attractive for coping with stressed operation conditions. One important advantage of NLP for OPF is that it naturally captures nonlinearities in stressed power system behavior, rendering them into an excellent tool for modelling and simulation of modern power system operations. As a result there is a great demand for new algorithms and software tools able to address strong nonlinearities in system behavior, in order to guarantee reliable and economic system operation. Additionally, in  order  to  comply  with  current  and  near  future  power systems  operating  conditions,  the  incorporation  of  flexible security criteria and the consideration of fast control actions due to growing renewables penetration should be the natural, very ambitious and very challenging targets for the first generation of real-time AC OPF under uncertainty tools.

\section{Security Constrained OPF}

The security constrained OPF (SCOPF)~\cite{SCOPFinvention}, is an extension of the OPF problem, which finds an optimal operational state but at the same time takes into account a set of security constraints arising from the operation of the system under a set of postulated contingencies. It guarantees the whole power system can work under the nominal long-term cost-efficient operation plan, but can also remain in the operational state when some of the predetermined contingencies occur.
However, each additional contingency corresponds to an additional set of constraints in the OPF problem specified for the associated power grid. The nominal scenario and all contingency states are coupled, rendering the whole problem computationally intractable for realistic size grids. 

The SCOPF optimization problem considered in this work is the ``preventive" SCOPF, although the same algorithmic improvements would apply also to the ``corrective" variant. SCOPF problem is formulated as:
\begin{subequations}
	\label{SCOPF}
	\begin{align}
	\label{SCOPF:objective}
	\minimize{\bftheta_{c},\bfv_{c},\bfp_{c},\bfq_{c}}  & \sum_{l=1}^{\Ng} a_l(\pes{l}{0})^2 + b_l\pes{l}{0} +  c_l \\ 
	\subject{\,\;}  & \forall c \in \{ c_0, c_1, \ldots, c_{N_c}\}: \nonumber \\
	\label{SCOPF:PowerFlow}
	& \PowerFlowSCOPF(\vc_c, \sgc_c) = \0, \\  
	\label{SCOPF:LineFlowLimits}
	& \LineFlowSCOPF(\vc_c)  \leq \Max{\F}_\text{L}, \\
	\label{SCOPF:VoltageLimits}
	& \Min{\bfv}      \leq \bfv_c      \leq \Max{\bfv},\quad \bftheta_c^{ref} = 0,  \\
	\label{SCOPF:PowerLimits}
	& \Min{\bfp} \leq \bfp_c  \leq \Max{\bfp},  \quad 
	\Min{\bfq} \leq \bfq_c  \leq \Max{\bfq}, \\
	\label{SCOPF:SCOPFVm}
	& \forall b \in \mathcal{B}_{PV} :\  \bfv_c = \bfv_{c_0} , \\
	\label{SCOPF:SCOPFPg}
	& \forall  g \in \mathcal{B}_{PV} :\  \bfp_c = \bfp_{c_0}. 
	\end{align}
\end{subequations}
Note that the SCOPF problem replicates the OPF constrains and variables for each contingency scenario $c$. The values of the non-automatic control variables are the same in all system scenarios, as expressed by the two non-anticipatory constraints \eqref{SCOPF:SCOPFVm} and \eqref{SCOPF:SCOPFPg}. These declare that the voltage magnitude and the active power generation at the PV buses $\mathcal{B}_{PV}$ (also known as generator buses, the active power and voltage magnitude are specified, therefore the name PV bus) should remain the same as in the nominal scenario $c_0$, regardless which contingency they are associated with. The only generator that is allowed to change its active power output is the generator at the reference bus $\mathcal{B}_{ref}$ (also known as slack or swing bus) as its active power generation can be modified to refill the power transmission losses occurring in each contingency $c$. This implies that part of the optimization vector $\bfx$ will be shared between the scenarios and part of it will be local to each contingency. Therefore, the vector of variables can be partitioned into local components $\bfx_c$ for each contingency $\forall c \in \{  c_0, c_1, \ldots, c_{N_c}\}$ and the global (shared) part $\bfx_g$:
\begin{align}
    \bfx_c &= [ \bftheta^\transpose, \bfv_i^\transpose, \bfq^\transpose, \bfp_j^\transpose]^\transpose,\ i \not \in \mathcal{B}_{PV}, j \in \mathcal{B}_{ref}, \label{eq:SCOPFxi}\\
   	\bfx_g &= [\bfv_i^\transpose,\bfp_j^\transpose ]^\transpose,\ i \in \mathcal{B}_{PV}, j \not \in \mathcal{B}_{ref}, \label{eq:SCOPFxg} \\
    \bfx &= [\bfx_0^\transpose, \bfx_1^\transpose,\ldots, \bfx_{N_c}^\transpose, \bfx_g^\transpose ]^\transpose. \label{eq:SCOPFx}
\end{align}
This ordering of the variables allow implementation of the efficient structure exploiting algorithms described in the following sections.

An important question arises in the SCOPF problem studies, which is how to select a reduced number of equivalent credible umbrella contingencies and associated variables that the SCOPF problem has to consider in order to obtain the same or nearly the same solution as with the full set of contingencies. In other words, how to identify contingencies that do not pose any security risks and at the same time do not restrict the optimal solution, therefore do not have to be considered in the SCOPF problem. The set of umbrella contingencies for a given SCOPF is strongly dependent on the operating conditions of the power grid (loading, grid configuration, etc.). Hence, as the parameters of the problem change, the membership in the set of umbrella contingencies also varies. Machine learning techniques are usually used for this purpose, including three classes of learning methods, namely machine learning, artificial neural networks and statistical pattern recognition \cite{ANN, Umbrella}. However, the resulting set might be still very large and sophisticated optimization algorithms are needed.
For the purpose of this study, such contingencies are selected that result in non-empty feasible region of the SCOPF problem \eqref{SCOPF}. This restricts the selection of the transmission lines that can experience failure and do not compromise operation of the power grid. Such lines are characterized by the following properties:
(i) no islands and isolated buses appear in the grid after the line failure, (ii) the reduced grid remains feasible in the PF sense, and (iii) only limited reactive power generation violations are allowed after the contingency occurrence. These considerations are usually made in the planning stage of the transmission grid design, but since synthetic power grid networks are used in this study, these considerations might not be satisfied for all conceivable contingencies.

Several techniques have been proposed in the literature aiming to reduce the computational complexity of the SCOPF problem, such as direct current (DC) approximations~\cite{SCOPFclustering}, filtering of contingencies using prior knowledge about the grid~\cite{SCOPFfiltering},~\cite{SCOPFfiltering1} or nonlinear decomposition of the grid into zones with the aid of additional user information~\cite{SCOPFdistributedUC}. Other algorithms adopt the dual decomposition or the  alternating  direction  method  of multipliers, e.g., applied to the DC OPF problems in~\cite{SCOPFadmm},\cite{decentralizedOPFGabyNEW}. Recently, in addition to classic CPU-based computing, graphical processing units (GPUs) have also been explored for security analysis in \cite{SSAgpu} or \cite{SCPFgpu}.
Alternatively, the SCOPF problem can be decomposed on the linear level. Interior point methods (IPMs) have been the most robust and successful tools for large-scale nonconvex optimization, and IPMs can easily exploit the problem structure~\cite{Gondzio2009}. Most of the computational time in IPMs is spent in the solution of the linear system, arising from the linearization of the Karush-Kuhn-Tucker (KKT) optimality conditions, the size of which grows with the number of considered contingencies and the size of the power grid. The KKT system is usually very large but sparse. Direct sparse linear solvers are the standard choice employed for the solution of the KKT. The solution can be accelerated using techniques tailored to the particular structure of the associated KKT linear system.

The parallel distributed Schur complement framework was applied to AC SCOPF problems and associated KKT systems in~\cite{Kang15}. The work focuses on  solving the Schur complement equations implicitly using a quasi-Newton preconditioned conjugate gradient method. However, the study does not evaluate the algorithm on large-scale problems and demonstrate scaling only up to 16 MPI processes. 
Structured non-convex optimization of large-scale energy systems using PIPS-NLP was performed in \cite{CosminMPSCOPF, CosminCurvature}. The parallel interior point optimization solver for nonlinear programming leverages the dual-block angular structure specific to the problem formulation by applying the Schur complement for efficient parallelization of the linear solves. It was illustrated how different model structures arise in power system domains and how these can be exploited to achieve high computational efficiency. Stochastic optimization problems on high-performance computers have been treated similarly in~\cite{olafAugmented, pardiso}. The numerical experiments suggest that supercomputers can be efficiently used to solve power grid optimization problems with thousands of scenarios under the strict time requirements of power grid operators, particularly due to improved linear algebra on a shared memory level.
The proposed additional Schur complement scheme in this work differentiates the proposed formulation from the one suggested by the previous work.  The benefits come from a more efficient direct sparse approach for the solution of the underlying sparse KKT system solved at each iteration of the IPM. The proposed method employs one additional Schur complement carefully chosen so that the sparsity of the KKT matrix is maintained. This way the size of the linear system decreases and this allows for memory savings and increased computational performance.

\section{Reduced Space IPM and Adjoint Formulation \label{chapter:adjoint}}

The optimization problems in the PDE-constrained problems (e.g., maximization of the oil production \cite{DrososConstraints}) are often solved by reduced space optimization methods \cite{Ghattas-01, Ghattas-02}. The discretized state variables $\bfx$ and the PDEs representing the nonlinear equality constraints \eqref{eq:PDEeq} are removed from the optimization problem and are treated explicitly during the evaluation of the objective function value and its gradient. Given the control variables $\bfu$, the system of the removed equality constraints is used to solve for the state variables $\bfx(\bfu)$. Thus, the equality constraints are implicitly satisfied. The efficient evaluation of the gradient information is achieved using the adjoint method. The computational cost of evaluating the inequality constraints gradients \eqref{eq:PDEineq} can be reduced using the constraint lumping techniques \cite{DrososConstraints}. The second order derivatives are usually not evaluated exactly due to the excessive computational cost, only approximations such as BFGS are used. The concepts of the reduced space IPM are first applied to the OPF problem, the extension to the SCOPF problems is discussed later.

The OPF problem
\begin{subequations}
	\label{eq:PDEoptim}
	\begin{align}
	\min_{\bfx(\bfu),\bfu} \qquad & f(\bfx(\bfu), \bfu)  \label{eq:PDEobj}\\
	\text{subject to} \qquad &  \bfg(\bfx(\bfu), \bfu) = \0, \label{eq:PDEeq}\\
	 &  \bfh(\bfx(\bfu), \bfu) \leq \0, \label{eq:PDEineq} \\
	 & \bfx^{min} \leq \bfx(\bfu) \leq \bfx^{max}, \label{eq:PDEstate}\\
	 & \bfu^{min} \leq \bfu \leq \bfu^{max}. \label{eq:PDEcontrol}
	\end{align}
\end{subequations}
can be interpreted in the context of PDE-constrained optimization problems, where the PF equations \eqref{eq:PDEeq} represent the role of the PDEs. The control variables $\bfu$ and the state variables $\bfx(\bfu)$ correspond to the known and unknown quantities in the PF problem. The known global control variables $\bfu$
contain real power produced at the generators (except the slack generator) and the bus voltage magnitude involved in reactive power balance. The local variables $\bfx(\bfu)$, representing the state of the power grid, 
are the remaining quantities such as voltage angles, voltage magnitudes at PQ buses, generator reactive power injections, etc.


The main advantage of the reduced space IPM is that the computational complexity mainly depends on the degrees of freedom, that is, the number of control variables $n_u=|\bfu|$. The reduced space technique is very suitable for solving the optimization problems where the number of control variables $n_x=|\bfx|$ is much smaller that the total number of variables, which is the case especially for the SCOPF problem. Furthermore, the reduced space approach is suitable for applications such as time critical optimization, where one may want to terminate an algorithm before optimality has been reached and be assured that the current approximate solution is feasible. In other applications, the objective function might not be defined outside of the feasible region. 
The main drawback is that in the line search framework, the objective function might be required to be evaluated multiple times with the updated control variables, which requires to resolve the dependent state variables repeatedly. In the OPF and SCOPF problems, such computation is expensive but the evaluation of the gradients will require much more computational resources.

In  the  reduced space  approach,  the  optimizer controls only $\bfu$ variables but since $f(\bfx(\bfu),\bfu)$ and $\bfh(\bfx(\bfu),\bfu)$ must be evaluated, the variables $\bfx(\bfu)$ must be known. The states $\bfx(\bfu)$ of the problem are obtained as the solution of \eqref{eq:PDEeq} for given controls $\bfu$. Therefore, the network constraints are not forwarded to the optimizer as additional equality constraints, which considerably reduces the size of the optimization problem, not only in terms of the primal state variables but also the associated dual variables for the constraints. The reduced space problem is formulated as
\begin{subequations}
	\label{eq:reducedPDEoptim}
	\begin{align}
	\min_{\bfu} \qquad & f(\bfx(\bfu), \bfu)  \label{eq:reducedPDEobj}\\
	\text{subject to} \qquad &  \bfh(\bfx(\bfu), \bfu) \leq \0, \label{eq:reducedPDEineq} \\
	 & \bfx^{min} \leq \bfx(\bfu) \leq \bfx^{max}, \label{eq:reducedPDEstate}\\
	 & \bfu^{min} \leq \bfu \leq \bfu^{max}. \label{eq:reducedPDEcontrol}
	\end{align}
\end{subequations}

\section{Evaluation of the State Variables \label{sec:adjointOrdering}}
The variables in the full space OPF problems are the complex voltages specified at each bus $\vc = \bfv e ^{j\bftheta}$ and the generator power injections $\sgc = \bfp + j \bfq$. We further distinguish three types of buses, namely reference REF, PV and PQ. The known control variables in the reduced space problem are $\bfu = (\bfv_{PV}, \bfp_{PV})$, that is, voltage magnitudes at the PV buses and the generator active powers at the PV buses. The unknown state variables $\bfx_1 = (\bftheta_{PV}, \bftheta_{PQ}, \bfv_{PQ})$ and $\bfx_2 = (p_{REF}, q_{REF}, \bfq_{PV})$ form the vector with the ordering $\bfx = (\bfx_1, \bfx_2)$. The remaining two variables, $\theta_{REF}, v_{REF}$, are fixed and known, therefore are not considered as additional variables.

The unknown state variables $\bfx$ are computed by the standard Newton iterations with a set of the nonlinear mismatch equations $\bfg(\bfx) = \0$, defined in \eqref{eq:PDEeq}, where the control variables $\bfu$ are considered as known, fixed parameters. The equations are ordered as $\bfg(\bfx) = (\bfg_1(\bfx), \bfg_2(\bfx))$, where $\bfg_1(\bfx) = (\Re(\bfg_{PV}), \Re(\bfg_{PQ}), \Im(\bfg_{PQ}))$ and $\bfg_2(\bfx) = (\Re(g_{REF}), \Im(g_{REF}), \Im(\bfg_{PV}))$.
Considering the aforementioned ordering of the state variables and the mismatch equations, the Newton iterations
\begin{equation}
    \frac{\partial \bfg(\bfx)}{\partial \bfx} \Delta \bfx  = - \bfg(\bfx), \label{eq:PFproblem}
\end{equation}
have the following structure
\begin{equation}
    \begin{bmatrix}
    \bfg_{11} & \bfg_{12} \\
    \bfg_{21} & \bfg_{22}
    \end{bmatrix}
    \begin{bmatrix}
    \Delta \bfx_1 \\
    \Delta \bfx_2
    \end{bmatrix}
     = -
    \begin{bmatrix}
    \bfg_1(\bfx) \\
    \bfg_2(\bfx)
    \end{bmatrix}.
\end{equation}
The block $\bfg_{12} = \partial \bfg_{1}(\bfx) / \partial \bfx_2$ is zero block, since none of the components in $\bfx_2$ are present in the equations $\bfg_1$. The  update $\Delta \bfx_1$ can be computed by using the standard power flow problem \cite{meier2006book}. If the power grid is connected, the power flow problem $\bfg_{11} \Delta \bfx_1=-\bfg_1(\bfx)$ is well posed, and the Jacobian $\bfg_{11} = \partial \bfg_{1}(\bfx) / \partial \bfx_1$ is square and non-singular for fixed $\bfu$. The update $\Delta \bfx_2$ can be computed from $\bfg_{22} \Delta \bfx_2= -\bfg_2(\bfx) - \bfg_{21}\Delta \bfx_1$, where $\bfg_{22} = -\bfI$ is negative identity matrix.

\section{Evaluation of the Gradients}
Usage of the gradient-based methods calls for an efficient tool to compute the gradient of the objective function with respect to the design variables, with maximum accuracy and minimum CPU cost. Using the chain rule, the  gradient $\frac{df}{d{\bfu}}$ can be expressed as
\begin{equation}
    \frac{df}{d{\bfu}} = \partial_{\bfu}f + \partial_{\bfx}f d_{\bfu}\bfx. \label{eq:chainF}
\end{equation}
The term $d_{\bfu}\bfx$ is not known, since $\bfx$ cannot be expressed analytically as a function of $\bfu$. The solution of \eqref{eq:PDEeq} can be obtained using only numerical methods. 
One method to approximate the gradient is to use $n_u$ finite differences over the elements of $\bfu$ in order to approximate $\frac{df}{d{\bfu}}$. This is computationally expensive, since it requires solving \eqref{eq:PDEeq} for each finite difference. 
The adjoint method addresses the problem of computing $\frac{df}{d{\bfu}}$ with the cost of one adjoint system solution, irrespective of the number of control variables $n_u$ in case of the OPF problem.

Many optimization methods make use also of the second order information, that is, not only the objective gradient but also Hessian is considered. The second order algorithms have a much faster theoretical convergence rate, but this comes with the cost of evaluating the Hessian. An overview of the Hessian evaluation cost is discussed in \cite{adjointHessian}. Alternatively, one may rely on first order Hessian approximations, e.g. BFGS of DFP \cite{NocedalBookOptim}.




\subsection{Derivation of the Adjoint Equations}
In order to use the gradient based optimization techniques to solve \eqref{eq:PDEoptim}, we need to know the first variation of the objective function $\delta f$.
The efficient computation is addressed by the adjoint method. It exploits the fact that during the solution process of the state variables \eqref{eq:PFproblem} the Jacobian matrix $\partial_x \bfg$ is used in Newton's method. The adjoint method uses the transpose of the Jacobian in order to compute the objective function's first variation.

We can introduce the augmented objective function
\begin{equation}
    \mathcal{L}(\bfx,\bfu,\bflambda) \equiv f(\bfx, \bfu) + \bflambda^{\transpose} \bfg(\bfx, \bfu),
    \label{eq:adjointObj}
\end{equation}
where $\bflambda$ is the vector of Lagrange multipliers. As $\bfg(\bfx, \bfu)$ is zero by construction, the multipliers may be chosen freely without changing the minimum of the augmented function, therefore, $f(\bfx, \bfu) = \mathcal{L}(\bfx,\bfu,\bflambda)$. Since $\bfg(\bfx,\bfu)$ is zero everywhere, its gradient is zero as well,
\begin{equation}
    d_u \bfg = \frac{\partial \bfg}{\partial \bfu} \delta \bfu + \frac{\partial \bfg}{\partial \bfx} \delta \bfx = \0,
\end{equation}
and the first variation of $f$ can be expressed as
\begin{equation}
   \delta f(\bfx,\bfu) = \delta \mathcal{L}(\bfx,\bfu,\bflambda) = \left( \frac{\partial f}{\partial \bfx} + \bflambda^{\transpose} \frac{\partial \bfg}{\partial \bfx} \right) \delta \bfx + \left(\frac{\partial f}{\partial \bfu} + \bflambda^{\transpose} \frac{\partial \bfg}{\partial \bfu} \right) \delta \bfu. \label{eq:d_adjointObj}
\end{equation}
In order to achieve $\delta \mathcal{L} = \0$ we require $\partial \mathcal{L}/\partial \bfx = \0$ and $\partial \mathcal{L}/\partial \bfu = \0$. To satisfy $\partial \mathcal{L}/\partial \bfx = \0$, we require that the Lagrange multipliers satisfy the following equations:
\begin{equation}
    \left(\frac{\partial \bfg}{\partial \bfx}\right)^{\transpose} \bflambda = -\left(\frac{\partial f}{\partial \bfx}\right)^{\transpose}.
    \label{eq:ls_adjointObj}
\end{equation}
With this choice of the Lagrange multipliers the first term of \eqref{eq:d_adjointObj} becomes zero and the gradient of the objective function with respect to the controls that is required by the optimization software is
\begin{equation}
    \frac{\delta f}{\delta \bfu} = 
    \frac{\partial f}{\partial \bfu} +  \bflambda^{\transpose} \frac{\partial \bfg}{\partial \bfu}.
    \label{eq:adjointObjReduced}
\end{equation}
Note that in the process of solving for the state variables \eqref{eq:PFproblem} the Jacobian of the constraints $\partial \bfg/\partial \bfx$ is used to determine the Newton step
\begin{equation}
    \label{eq:adjointNewton}
    \frac{\partial \bfg}{\partial \bfx} \Delta x  = - \bfg.
\end{equation}
The adjoint equation \eqref{eq:ls_adjointObj} solves a linear system that differs in form from \eqref{eq:adjointNewton} only by the adjoint operation and the RHS vector. If the power network is connected, the power flow problem $\bfg(\bfx,\bfu)=\0$ is square and well posed, and the Jacobian is nonsingular for fixed $\bfu$.

\subsection{Detailed Treatment of the Objective Adjoint System \label{sec:detailedAdjoint}}

We want to evaluate the objective function gradient \eqref{eq:d_adjointObj}, where the ordering of the variables and the constraints is defined in section \ref{sec:adjointOrdering}. The gradient, considering also the contribution of the augmented function, is
\begin{equation}
   \delta f(\bfx,\bfu) = \delta \mathcal{L}(\bfx,\bfu,\bflambda) = \left( \frac{\partial f}{\partial \bfx} + \bflambda^{\transpose} \frac{\partial \bfg}{\partial \bfx} \right) \delta \bfx + \left(\frac{\partial f}{\partial \bfu} + \bflambda^{\transpose} \frac{\partial \bfg}{\partial \bfu} \right) \delta \bfu. \nonumber
\end{equation}
Since the objective function depends only on the $\bfp$ variables, there is only a single nonzero entry in the partial derivative w.r.t. the state variables $\bfx$, and a zero and a dense block in partial derivative w.r.t the control variables
\begin{equation}
\frac{\partial f}{\partial \bfx} =
\begin{bmatrix}
  \frac{\partial f}{\partial \bfx_1} \\
  \frac{\partial f}{\partial \bfx_2}
\end{bmatrix},
\qquad
\frac{\partial f}{\partial \bfu} =
\begin{bmatrix}
  \frac{\partial f}{\partial \bfv_{PV}} \\
  \frac{\partial f}{\partial \bfp_{PV}}
\end{bmatrix},
\end{equation}
where $\partial f/\partial \bfx_1 = \0$ and $\partial f/\partial \bfx_2 = [nnz, 0, 0, \ldots]$ where the nonzero element corresponds to $\partial f/\partial p_{REF}$. Similarly, $\partial f/\partial \bfv_{PV} = \0$ and $\partial f/\partial \bfp_{PV} \neq \0$ is a dense block of nonzeros.

\subsubsection*{Evaluation of the Multipliers}
In order to evaluate the objective function gradient \eqref{eq:d_adjointObj}, we require that the term in front of $\delta \bfx$ evaluates to zero by an appropriate choice of the multipliers $\bflambda = [\bflambda_1, \bflambda_2]$. We require that
\begin{equation}
    \frac{\partial f}{\partial \bfx} + \bflambda^{\transpose} \frac{\partial \bfg}{\partial \bfx} = 0,
\end{equation}
thus obtaining the adjoint system,
\begin{equation}
    \frac{\partial \bfg}{\partial \bfx}^\transpose \bflambda = -\frac{\partial f}{\partial \bfx}.
\end{equation}
The structure of the adjoint linear system $(\partial \bfg / \partial \bfx)^\transpose$, considering ordering of the variables and equations from section \ref{sec:adjointOrdering}, and the fact that $\bfg_{12} = \0$, $\bfg_{22} = -\bfI$, is the following
\begin{equation}
    \begin{bmatrix}
    \bfg_{11}^\transpose & \bfg_{21}^\transpose \\
    \0 & -\bfI
    \end{bmatrix}
    \begin{bmatrix}
    \bflambda_1 \\
    \bflambda_2
    \end{bmatrix}
     = -
    \begin{bmatrix}
    \frac{\partial f}{\partial \bfx_1} \\
    \frac{\partial f}{\partial \bfx_2}
    \end{bmatrix}.
\end{equation}
We can decouple the solution of $\bflambda_1$ and $\bflambda_2$ by solving first for $\bflambda_2$, which is trivial. Afterwards, the smaller adjoint system is solved
\begin{align}
    \bflambda_2 &= \frac{\partial f}{\partial \bfx_2}, \\
    \bfg_{11}^\transpose \bflambda_1 &= -\underbrace{\frac{\partial f}{\partial \bfx_1}}_{=\0} - \bfg_{21}^\transpose \underbrace{\bflambda_2}_{\partial f/\partial \bfx_2}.
\end{align}
Furthermore, by considering the structure of the RHS vector $\partial f/\partial \bfx_2$, which has only a single nonzero entry at the first position, only the first column of $\bfg_{21}^\transpose$ needs to be evaluated, which corresponds to the gradient of a single constraint $\partial \Re(g_{REF})/\partial \bfx_1$.

\subsubsection*{Evaluation of the Gradient}
Once we have the value of the multipliers, we substitute the multipliers $\bflambda$ in the second term in \eqref{eq:d_adjointObj} in order to evaluate the gradient, that is
\begin{equation}
\delta f(\bfx,\bfu) = \delta \mathcal{L}(\bfx,\bfu,\bflambda) = \frac{\partial f}{\partial \bfu} + \bflambda^{\transpose} \frac{\partial \bfg}{\partial \bfu}.    
\end{equation}
The nonzero pattern of the multipliers is $\bflambda_2 = [nnz, 0, 0, \ldots]$, therefore we need to evaluate only the rows of the $\partial \bfg/\partial \bfu$ corresponding to gradients of constraints $[ \bfg_1, \Re(g_{REF}) ]$. Furthermore, $\partial \bfg / \partial \bfp_{PV}$ is zero, except $\partial \Re(\bfg_{PV}) / \partial \bfp_{PV} = -\bfI$.

\subsection{Constraint Handling}
Constraints that appear as simple bound constraints on the control variables \eqref{eq:PDEcontrol}, such as generator power output bounds, can be used directly as inputs to the optimizer. These constraints are removed internally by the optimizer and transformed into standard logarithmic barrier terms. Constraints on state variables \eqref{eq:PDEstate}, such as the reactive power output, require the solution of the forward PF problem for the evaluation of the constraint and the solution of the adjoint problem for the evaluation of the gradient of the constraint with respect to the control variables. 

The number of nonlinear inequality constraints \eqref{eq:PDEineq} specified in the standard OPF control problem are $2\Nl$, where $\Nl$ is the number of transmission lines. The inequality constraints consist of two constraints for each line, one for the ``from" end and one for the ``to" end, imposing limits on apparent power flows expressed in MVA. Computing the gradient for each one of these constraints requires the evaluation of the Lagrange multipliers corresponding to each constraint. Similarly to the equation \eqref{eq:adjointObj}, the augmented constraint $\mathcal{H}^i$ is formulated for each of the inequality constraints, $i=1 \ldots 2\Nl$,
    \begin{equation}
       \delta \mathcal{H}^i(\bfx,\bfu,\bflambda) = \left( \frac{\partial h^i}{\partial \bfx} + \bflambda^{\transpose} \frac{\partial \bfg}{\partial \bfx} \right) \delta \bfx + \left(\frac{\partial h^i}{\partial \bfu} + \bflambda^{\transpose} \frac{\partial \bfg}{\partial \bfu} \right) \delta \bfu. \label{eq:d_adjointIneq}
    \end{equation}
The first term is required to be zero, therefore the same adjoint system is solved for $\bflambda$ as in \eqref{eq:d_adjointObj}, but with a different RHS. Since \eqref{eq:d_adjointIneq} requires the solution of a linear system, $2\Nl$ linear systems have to be solved in total to evaluate the gradient of the nonlinear constraints.

The bounded state variables are voltage magnitudes at PQ and reference buses, active power at the reference generator, and the reactive power produced at generators. The gradient evaluation of the state variable bounds \eqref{eq:PDEstate} is treated similarly as the nonlinear inequality constraints. We define the upper bound constraints as $\bfc_{max}(\bfx) = \bfx - \bfx^{max}$ and the first variation of the augmented constraint is
\begin{equation}
   \delta \mathcal{C}^i_{max}(\bfx,\bfu,\bflambda) = \left( \frac{\partial c_{max}^i}{\partial \bfx} + \bflambda^{\transpose} \frac{\partial \bfg}{\partial \bfx} \right) \delta \bfx + \left(\underbrace{\frac{\partial c_{max}^i}{\partial \bfu}}_{=\0} + \bflambda^{\transpose} \frac{\partial \bfg}{\partial \bfu} \right) \delta \bfu. \label{eq:d_adjointBounds}
\end{equation}
Since ${\partial \bfc_{max}}/{\partial \bfx}$ is $-\bfI$, the RHS vectors of the adjoint systems will be columns of the identity matrix. The lower bound constraints $\bfc_{min}(\bfx) = \bfx^{min} - \bfx$ are treated in the same way. The total number of adjoint systems solved is $2n_x$, however, since the lower bound adjoint system varies solely in the sign of the RHS, only $n_x$ systems need to be solved.

\subsection{Detailed Treatment of the Constraints Adjoint System}
In order to evaluate the constraint adjoint \eqref{eq:d_adjointIneq}, we follow similar process as in case of the objective function. That is, we require the term in front of $\delta \bfx$ to be zero
\begin{equation}
   \delta \mathcal{H}^i(\bfx,\bfu,\bflambda) = \left( \frac{\partial h^i}{\partial \bfx} + \bflambda^{\transpose} \frac{\partial \bfg}{\partial \bfx} \right) \delta \bfx + \left(\frac{\partial h^i}{\partial \bfu} + \bflambda^{\transpose} \frac{\partial \bfg}{\partial \bfu} \right) \delta \bfu. \nonumber
\end{equation}
We need to solve the adjoint system, similarly as in the section \ref{sec:detailedAdjoint}. The adjoint system is identical, only the RHS vector vary.

The nonlinear constraints $h^i(\bfx, \bfu), i=1\ldots,\Nl$ are functions of the $\bftheta$ and $\bfv$, thus the nonzero entries in the partial derivative w.r.t. the state variables $\bfx$, and the control variables $\bfu$ are
\begin{equation}
\frac{\partial h^i}{\partial \bfx} =
\begin{bmatrix}
  \frac{\partial h^i}{\partial \bfx_1} \\
  \frac{\partial h^i}{\partial \bfx_2}
\end{bmatrix},
\qquad
\frac{\partial h^i}{\partial \bfu} =
\begin{bmatrix}
  \frac{\partial h^i}{\partial \bfv_{PV}} \\
  \frac{\partial h^i}{\partial \bfp_{PV}}
\end{bmatrix},
\end{equation}
where $\partial h^i/\partial \bfx_1$ contains nonzero entries and $\partial h^i/\partial \bfx_2 = \0$. Similarly, $\partial h^i/\partial \bfv_{PV}$ may or may not contain nonzero entries, which is given by the power grid connectivity  (e.g. if the PV bus is connected to one, or both ends of the $i$th transmission line). Finally, $\partial f/\partial \bfp_{PV} = \0$.

Assuming this nonzero structure, the adjoint system
\begin{equation}
    \begin{bmatrix}
    \bfg_{11}^\transpose & \bfg_{21}^\transpose \\
    \0 & -\bfI
    \end{bmatrix}
    \begin{bmatrix}
    \bflambda_1 \\
    \bflambda_2
    \end{bmatrix}
     = -
    \begin{bmatrix}
    \frac{\partial h^i}{\partial \bfx_1} \\
    \frac{\partial h^i}{\partial \bfx_2}
    \end{bmatrix},
\end{equation}
simplifies to 
\begin{align}
    \bflambda_2 &= \0, \\
    \bfg_{11}^\transpose \bflambda_1 &= -\frac{\partial h^i}{\partial \bfx_1}.
\end{align}
The adjoint gradient can be subsequently evaluated substituting the multipliers $\bflambda$ into 
\begin{equation}
\delta h^i(\bfx,\bfu) = \delta \mathcal{H}^i(\bfx,\bfu,\bflambda) = \frac{\partial h^i}{\partial \bfu} + \bflambda^{\transpose} \frac{\partial \bfg}{\partial \bfu}.    
\end{equation}
The nonzero pattern of the multipliers, specifically $\bflambda_2 = \0$, allows us to evaluate only the rows of the $\partial \bfg/\partial \bfu$ corresponding to gradients of constraints $\bfg_1$ and ignore all gradients of constraint $\bfg_2$. Furthermore, $\partial \bfg / \partial \bfp_{PV}$ is zero, except $\partial \Re(\bfg_{PV}) / \partial \bfp_{PV} = -\bfI$.
\section{Constraint Lumping}
Explicit evaluation of all inequality constraints \eqref{eq:PDEineq} and especially their gradients is computationally intensive for large, realistically sized problems. It suffices for the feasible solution that the maximum of all lines with given rating, $l=1\ldots 2\Nl$, is within the limits, that is,
\begin{equation}
 \max_{l}(h^l) \leq 0. \label{eq:lumpingLines}
\end{equation}
If this maximum is honored, than all of the true constraints are guaranteed to be satisfied. Constraints that are described by nondifferentiable functions can be challenging to incorporate in the optimization frameworks. The approximation of the max function by smooth functions is known in the literature as constraint lumping \cite{ConorGas,DrososConstraints}. The max function from \eqref{eq:lumpingLines} is approximated as
\begin{equation}
    c = \max_{l}(h^l) \approx \alpha \ln{Q}, \label{eq:lumpingApproxLines}
\end{equation}
where $\alpha = 0.05h^{max}$ ({which needs to be chosen once at the beginning and fixed, not updated in every iteration with new max(h), additionally it needs to be positive}) and
\begin{equation}
    Q = \sum_{l=1}^{2\Nl} e^{h^l/\alpha}.
\end{equation}

All the line constraints \eqref{eq:PDEineq} in the OPF problem can be thus replaced by a single augmented constraint \eqref{eq:lumpingApproxLines} (assuming that all lines have the same power flow rating, otherwise there will be as many max constraints as there are line categories). The gradient of the constraint, which is required by the optimizer, augmented with the equality constraints \eqref{eq:PDEeq}, is given by
\begin{align}
\partial (c + \bflambda^{\transpose} \bfg)
    &= \left( \frac{\alpha}{Q} \frac{\partial Q}{\partial \bfu} + \bflambda^{\transpose} \frac{\partial \bfg}{\partial \bfu} \right) \partial \bfu + \left( \frac{\alpha}{Q} \frac{\partial Q}{\partial \bfx} + \bflambda^{\transpose} \frac{\partial \bfg}{\partial \bfx} \right) \partial \bfx \nonumber \\
    &= \left( \frac{\alpha}{Q} \sum_{i=1}^{2\Nl} e^{h^l/\alpha} \frac{\partial h^l/\alpha}{\partial \bfu} + \bflambda^{\transpose} \frac{\partial \bfg}{\partial \bfu} \right) \partial \bfu + \left( \frac{\alpha}{Q} \sum_{i=1}^{2\Nl} e^{h^l/\alpha} \frac{\partial h^l/\alpha}{\partial \bfx} + \bflambda^{\transpose} \frac{\partial \bfg}{\partial \bfx} \right) \partial \bfx \nonumber \\
    &= \left( \frac{1}{Q} \sum_{i=1}^{2\Nl} e^{h^l/\alpha} \frac{\partial h^l}{\partial \bfu} + \bflambda^{\transpose} \frac{\partial \bfg}{\partial \bfu} \right) \partial \bfu + \left( \frac{1}{Q} \sum_{i=1}^{2\Nl} e^{h^l/\alpha} \frac{\partial h^l}{\partial \bfx} + \bflambda^{\transpose} \frac{\partial \bfg}{\partial \bfx} \right) \partial \bfx.
    \label{eq:adjoingLumpedLines}
\end{align}
The second term in \eqref{eq:adjoingLumpedLines} is required to be zero, therefore, the Lagrange multipliers $\bflambda$ for the constraints are computed from the solution of the  adjoint problem
\begin{equation}
\frac{\partial \bfg^{\transpose}}{\partial \bfx} \bflambda =   -\frac{1}{Q} \sum_{i=1}^{2\Nl} e^{h^l/\alpha} \frac{\partial h^l}{\partial \bfx}.
\end{equation}
Note that only a single adjoint system needs to be solved in order to evaluate the gradient of the lumped line flow constraints, compared to $2\Nl$ adjoint systems in \eqref{eq:d_adjointIneq}. It is important to recognize that the approach used for constraint lumping can impact the performance of the optimization procedure. The smaller the coefficient multiplying $h^{max}$ in \eqref{eq:lumpingApproxLines}, the more accurate the approximation of max becomes. It is important to use a small numerical value, e.g. $0.05$, in order to ensure that no overflow occurs in any of the exponential terms in the summation. The approximation of max in \eqref{eq:lumpingApproxLines} is always greater than the maximum of the component line flows, so if this maximum is honored, the true constraint is guaranteed to be satisfied.

\section{Extension to Security Constrained OPF}
The formulation of the problem and computation of the adjoint systems in the SCOPF problem need to account for multiple sets of the state variables, corresponding to multiple contingency scenarios. The formulation of the problem \eqref{eq:PDEoptim} is extended by the additional sets of control variables $\bfx_c(\bfu)$ and constraints for each contingency scenario $c = 0,1,\ldots, NC$ including also the nominal case. The vector of all state variables is $\bfx = (\bfx_0^{\transpose}, \bfx_1^{\transpose}, \ldots, \bfx_{NC}^{\transpose})^{\transpose}$. The SCOPF problem is formulated as
\begin{subequations}
	\label{eq:PDEoptimSCOPF}
	\begin{align}
	\min_{\bfu} \qquad & f(\bfx(\bfu), \bfu)  \label{eq:PDEobjSCOPF}\\
	\text{subject to} \qquad & \forall c \in \{0,1,\ldots, NC\}: \nonumber \\
	 &  \bfh_c(\bfx_c(\bfu), \bfu) \leq \0, \label{eq:PDEineqSCOPF} \\
	 & \bfx^{min} \leq \bfx_c(\bfu) \leq \bfx^{max}, \label{eq:PDEstateSCOPF}\\
	 & \bfu^{min} \leq \bfu \leq \bfu^{max}. \label{eq:PDEcontrolSCOPF}
	\end{align}
\end{subequations}
The state variables $\bfx_c(\bfu)$ are computed by solving NC+1 PF problems \eqref{eq:PFproblem}. The set of the nonlinear equations $\bfg_c(\bfx_c) = \0$ is slightly different for each scenario, representing the modified grid with the particular contingency $c$. Note that the computations are independent of each other, therefore can be evaluated in parallel.

The augmented objective function used to formulate the adjoint system is
\begin{equation}
    \mathcal{L}(\bfx,\bfu,\bflambda_c) \equiv f(\bfx, \bfu) + \sum_{c} \bflambda_c^{\transpose} \bfg_c(\bfx_c, \bfu),
    \label{eq:adjointObjSCOPF}
\end{equation} and its first variation is
\begin{equation}
   \delta \mathcal{L}(\bfx,\bfu,\bflambda_c) = \sum_{c} \left( \frac{\partial f}{\partial \bfx_c} + \bflambda_c^{\transpose} \frac{\partial \bfg_c}{\partial \bfx_c} \right) \delta \bfx_c + \left(\frac{\partial f}{\partial \bfu} + \sum_{c}\bflambda_c^{\transpose} \frac{\partial \bfg_c}{\partial \bfu} \right) \delta \bfu. \label{eq:d_adjointObjSCOPF}
\end{equation}
Note the term ${\partial f}/{\partial \bfx_c}$ in the equation \eqref{eq:d_adjointObjSCOPF} above. Since the standard SCOPF objective function is defined as the cost of the active power generation, and the active power outputs are the control variables $\bfu$, the alluded term will be a vector of zeros for all scenarios. Only the active power output of the reference generator is part of the state variables. The reference generator needs to account for different transmission losses in every scenario. If we include the reference generator power output in the cost function only for the nominal case ($c=0$), then only a single adjoint system needs to be solved in order to evaluate the gradient of the objective function since the other terms vanish ($\frac{\partial f}{\partial \bfx_c}=\0$, therefore, $\bflambda_c = \0$ for $c>0$). If the reference generator active power output from every scenario is considered in the objective function, then $NC+1$ adjoint systems need to be solved in order to evaluate the gradient of the objective function. For comparison, there was one adjoint system \eqref{eq:d_adjointObj} in the OPF problem, which corresponds to considering all generators power output for the nominal case.

The augmented nonlinear line power flow constraint function for the line $i$ in contingency scenario $c$ is defined as \eqref{eq:adjointIneqSCOPF}. A similar process applies to the state variable bounds.
\begin{equation}
    \mathcal{H}^i_c(\bfx_c,\bfu,\bflambda_c) \equiv h_c^i(\bfx_c, \bfu) + \bflambda_c^{\transpose} \bfg_c(\bfx_c, \bfu).
    \label{eq:adjointIneqSCOPF}
\end{equation}
The gradient is expressed as
\begin{equation}
   \delta \mathcal{H}_c^i(\bfx_c,\bfu,\bflambda_c) = \left( \frac{\partial h_c^i}{\partial \bfx_c} + \bflambda_c^{\transpose} \frac{\partial \bfg_c}{\partial \bfx_c} \right) \delta \bfx_c + \left(\frac{\partial h_c^i}{\partial \bfu} + \bflambda_c^{\transpose} \frac{\partial \bfg_c}{\partial \bfu} \right) \delta \bfu. \label{eq:d_adjointIneqSCOPF}
\end{equation}

The computational complexity for the SCOPF problem grows as the problem size increases and as more contingency scenarios are considered. Overall, a single adjoint systems needs to be solved for each one of the $2\Nl$ line constraints in each contingency scenario $c$, including the nominal case. The gradient of state variable bounds can be evaluated with the cost of solving $n_x$ adjoint systems for each scenario. Overall, for the standard OPF, $1+2\Nl+n_x$ adjoint systems need to be solved; in case of  SCOPF there are $(NC+1)(1+2\Nl+n_x)$ adjoint systems that need to be solved. Introducing the line flow constraints \eqref{eq:PDEineqSCOPF} for individual contingency scenarios may lead to an excessively large number of constraints for realistic networks. This can have an adverse effect for the overall runtime performance of the SCOPF problem. The constraint lumping can be applied, as proposed in the previous section, which reduces the computational burden. The computational time can be further reduced by parallel computation of the gradient information of the constraints associated with each contingency scenario. 
\section{Numerical Results}

\begin{figure}[!th]
	\centering

\begin{tikzpicture}


\pgfplotstableread{
ID COMPONENT	    LUMPED	        FULLBFGS	FULLHESSIAN
0 {Grad Obj.}         	2.00E+01		2.00E-01	2.00E-01
1 {Jac. Constraints}	1.42E+02		2.09E+02	2.08E+02
2 {Hessian}	    0.00E+00		0.00E+00	1.02E+02
3 {Eval. State}	1.08E+03		0.00E+00	0.00E+00
4 {KKT Solve}	1.20E+01		5.77E+02	1.44E+03
5 {Overall} 	1.26E+03		7.93E+02	1.75E+03
}\datatable;

\begin{axis}[name=symb,
width=\columnwidth, height=4.5cm, 
ybar,
bar width=7pt, 
ymode=log, 
log origin=infty,
ymin=1e-1,
ymax=3e3,
axis lines*=left, 
ymajorgrids, yminorgrids,
xticklabels from table={\datatable}{COMPONENT},
xtick={0, 1, 2, 3, 4, 5, 6, 7, 8},
xticklabel style={align=center, rotate=20, xshift=-0.0cm, anchor=north, font=\scriptsize},
ytick={0.1,1,10,100,1000,10000,100000},
yticklabel style={font=\scriptsize},
ylabel style={font=\scriptsize},
xlabel style={font=\scriptsize},
ylabel={Time (ms)},
legend style={at={(1.0,1.12)},legend cell align=left,align=right,draw=none,font=\tiny,legend columns=5}
]
\addplot[fill=mycolor0] table[x=ID, y=LUMPED] {\datatable};
\addlegendentry{Reduced Space IPM}

\addplot[fill=mycolor1, postaction={pattern=north east lines}] table[x=ID, y=FULLBFGS] {\datatable}; 
\addlegendentry{IPM BFGS}

\addplot[fill=mycolor2, postaction={pattern=north west lines}] table[x=ID, y=FULLHESSIAN] {\datatable}; 
\addlegendentry{IPM Exact}
\end{axis}
\end{tikzpicture}
	\caption{Comparison of the time complexity of IPM components (ms) using case9214pegase. \label{fig:adjointComponents}}
\end{figure}
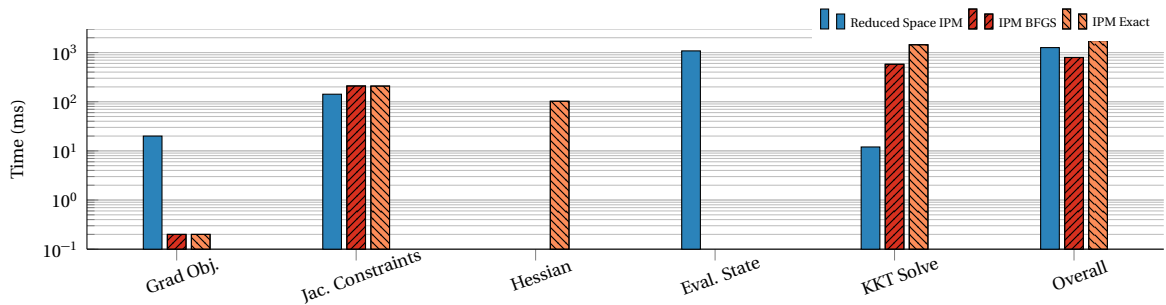

The most expensive component in the full space IPM with the exact Newton is the solution of the large-scale KKT system. While this step is avoided if using the BFGS approximation, the computational complexity of the reduced space quasi-Newton approach is expected to be dominated by evaluation of the state variables and evaluation of the gradients of the objective and constraint functions using the adjoint method, as demonstrated in Figure \ref{fig:adjointComponents}. %
The performance comparison for pegase9241 benchmark was performed for standard full space method with exact hessian, full space method with the BFGS Hessian approximation and the reduced space method with constraint lumping. The pegase9241 is the largest benchmark available, where the branch power flow limits are specified. The benchmark was run for 3 IPOPT iterations, where the number of line searches was always one, so that the number of function evaluation for all three methods was the same.
The precision of the timer is 1ms, so the routines taking less than one millisecond are not shown on the logarithmic plot. Similarly, the OPF Hessian is evaluated only for the exact IPM and the state variables evaluation is performed only in the reduced space method, therefore shown only for the corresponding methods.

Overall, the cost of the three approaches is very similar. The advantage of the reduced space approach is that the expensive components are embarrassingly parallel, thus the solution can greatly benefit from utilizing the parallel computations. However, by dropping the exact Hessian information, the convergence of the IPM might be negatively influenced.

\subsection*{Evaluation of the State Variables}
\begin{figure}[!th]
	\centering

\begin{tikzpicture}

\pgfplotstableread{
ID	TOL iters
0	5e-14 128
1	5e-13 135
2	5e-12 211
3	5e-11 261
}\datatable;

\begin{axis}[name=symb,
width=\columnwidth, height=4.5cm, 
ybar,
bar width=15pt, 
ymin=100,
ymax=300,
axis lines*=left, 
ymajorgrids, yminorgrids,
xticklabels from table={\datatable}{TOL},
xtick={0, 1, 2, 3, 4, 5, 6, 7, 8},
xticklabel style={align=center, rotate=20, xshift=-0.0cm, anchor=north, font=\scriptsize},
ytick={100,120, 140,...,300},
yticklabel style={font=\scriptsize},
ylabel style={font=\scriptsize},
xlabel style={font=\scriptsize},
ylabel={Iterations},
xlabel={PF tolerance},
legend style={at={(1.0,1.12)},legend cell align=left,align=right,draw=none,font=\small,legend columns=5}
]
\addplot[fill=mycolor0] table[x=ID, y=iters] {\datatable};
\addlegendentry{case118}
\end{axis}
\end{tikzpicture}
	\caption{IPM convergence sensitivity to the PF tolerance used for evaluation of the state variables. \label{fig:PFsensitivity}}
\end{figure}
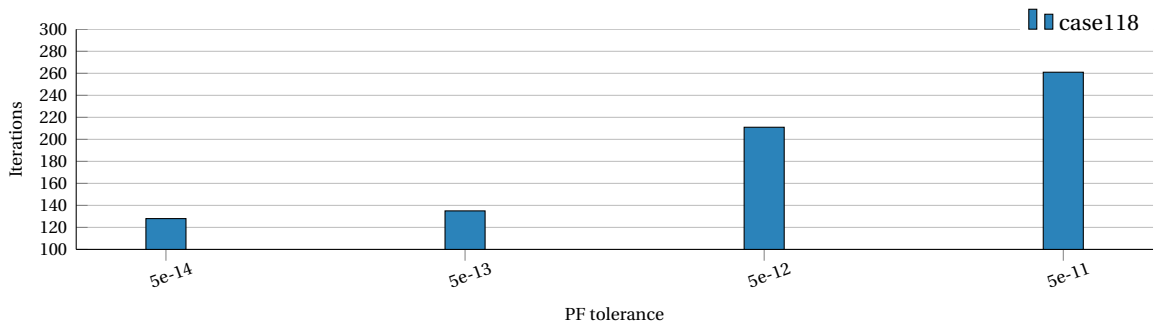

The evaluation of the state variables requires solving the power flow (PF) equations. The PF equations are backbone of the reduced space method. Newton-Raphson (NR) solution algorithm is used, due to its precision and fast convergence. It was observed that the reduced space method is sensitive to the quality of the PF solution, as demonstrated in Figure \ref{fig:PFsensitivity}. If the PF tolerance is relaxed, the IPM method needs significantly more iterations until convergence. This
possibly rules out computationally less expensive approximations such as Fast Decoupled Load-Flow (FDLF) method, which is based on NR, but greatly reduces its computational cost by means of a decoupling approximation that is valid in most transmission networks (the problem is split into two subproblems with smaller Jacobians which need to be factorized only once). The PF solution in case118 is solved up to the tolerance 1e-14, however for larger networks, such as PEGASE1354, only 5e-12 can be reached.

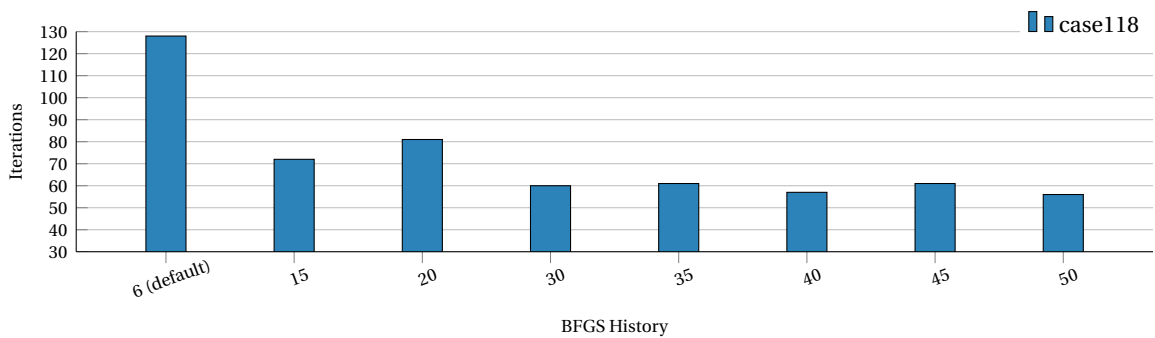
\begin{figure}[!th]
	\centering

\begin{tikzpicture}

\pgfplotstableread{
ID	History iters
0	{6 (default)} 128
1	15 72
2	20 81
3	30 60
4	35 61
5	40 57
6	45 61
7	50 56
}\datatable;

\begin{axis}[name=symb,
width=\columnwidth, height=4.5cm, 
ybar,
bar width=15pt, 
ymin=30,
ymax=130,
axis lines*=left, 
ymajorgrids, yminorgrids,
xticklabels from table={\datatable}{History},
xtick={0, 1, 2, 3, 4, 5, 6, 7, 8},
xticklabel style={align=center, rotate=20, xshift=-0.0cm, anchor=north, font=\scriptsize},
ytick={30,40,...,130},
yticklabel style={font=\scriptsize},
ylabel style={font=\scriptsize},
xlabel style={font=\scriptsize},
ylabel={Iterations},
xlabel={BFGS History},
legend style={at={(1.0,1.12)},legend cell align=left,align=right,draw=none,font=\small,legend columns=5}
]
\addplot[fill=mycolor0] table[x=ID, y=iters] {\datatable};
\addlegendentry{case118}
\end{axis}
\end{tikzpicture}
	\caption{IPM convergence sensitivity to the BFGS history length. \label{fig:Historysensitivity}}
\end{figure}

It was also observed that the reduced space method is sensitive to the quality of the limited memory BFGS approximation, as demonstrated in Figure \ref{fig:Historysensitivity}. If longer history of the problem is considered, the IPM method needs significantly less iterations until convergence.

\begin{figure}[!th]
	\centering

\begin{tikzpicture}

\pgfplotstableread{
ID	Barrier iters
0	{Adaptive} 128
1	{Monotone 1e-2} 500
2	{Monotone 5e-2} 481
3	{Monotone 9e-2} 311
4	{Monotone 1e-1} 342
5	{Monotone 5e-1} 500
6	{Monotone 1e-0} 500
}\datatable;

\begin{axis}[name=symb,
width=\columnwidth, height=4.5cm, 
ybar,
bar width=15pt, 
ymin=100,
ymax=550,
axis lines*=left, 
ymajorgrids, yminorgrids,
xticklabels from table={\datatable}{Barrier},
xtick={0, 1, 2, 3, 4, 5, 6, 7, 8},
xticklabel style={align=center, rotate=10, xshift=-0.0cm, anchor=north, font=\scriptsize},
ytick={100,150,...,500},
yticklabel style={font=\scriptsize},
ylabel style={font=\scriptsize},
xlabel style={font=\scriptsize},
ylabel={Iterations},
xlabel={Barrier Update Strategy},
legend style={at={(1.0,1.12)},legend cell align=left,align=right,draw=none,font=\small,legend columns=5}
]
\addplot[fill=mycolor0] table[x=ID, y=iters] {\datatable};
\addlegendentry{case118}
\end{axis}
\end{tikzpicture}
	\caption{IPM convergence sensitivity to the barrier parameter update strategy. \label{fig:Barriersensitivity}}
\end{figure}
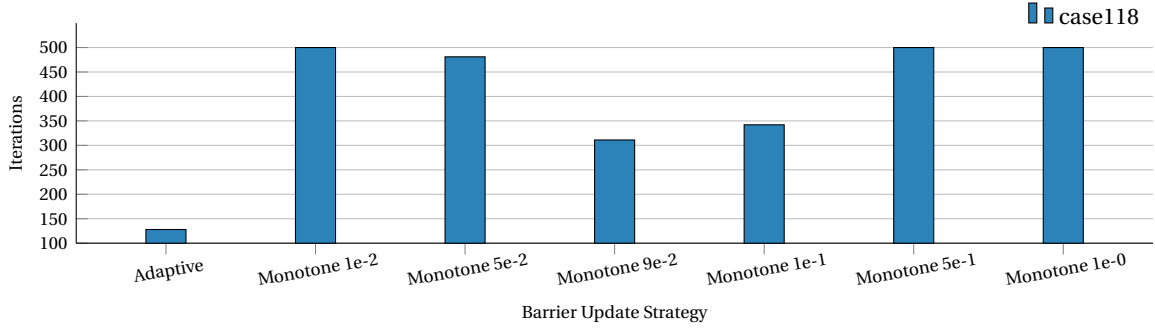

The reduced space IPM is also sensitive the the barrier parameter update strategy, as demonstrated in Figure \ref{fig:Barriersensitivity}. The monotone strategy results in significantly more iterations required until convergence.

\subsection*{Constraint Lumping}

The computational complexity of the reduced space IPM can be reduced by lumping the constraints representing the transmission lines with the same power flow limit. In the worst case, when each line has a different limit, the complexity would be improved by a factor of two, since there are two constraints for each line (i.e., the two constraints for each line will be replaced by a single constraint, approximating a maximum of the two). The research question is how does the approximation of the maximum influence the convergence of the reduced space IPM. 

It is important to note that the state variables (computed as a solution of the PF problem given the control variables) is very sensitive to the perturbation of the control variables. Even for a small perturbation of the controls, the output reactive powers $\bfq$ vary in great degree (while the max perturbation of controls by IPOPT is 2\%, the PF solution $\bfq$ change by a factor of 60) resulting in serious $\bfq$ bound violations and related branch power flow heat limits, representing the constraints that are subject to the lumping.

Due to the aforementioned sensitivity of the PF solution to perturbation of the control variables, the output powers $\bfq$ change rapidly. As a consequence, the branch power flows change similarly since they are function of the bus powers. Consequently, the values of the branch power flow violations across individual branches vary significantly, both across individual branches are between the iterations. Three variants of the nonlinear inequality constraints $\LineFlow(x)$  and their scaling were investigated, such that $\LineFlow(x)=SS^*-{\Max{\F}_\text{L}}^2 \leq 0$, where $S$ is the complex bus power and the branch power flow limit is denoted as $\Max{\F}_\text{L}$. The original constraints and two scaled versions are the following:
\begin{enumerate}
    \item $\LineFlow(x)$  : The maximum value is in order 6e4, minimum -2e2
    \item $\LineFlow(x)/{\Max{\F}_\text{L}}^2$: The maximum value is in order 7e2, minimum -1
	\item $\frac{\LineFlow(x)/({\Max{\F}_\text{L}}^2 )}{\max(\LineFlow(x)/({\Max{\F}_\text{L}}^2 ))}$: The maximum is always 1
\end{enumerate}
	
In cases 2. and 3. there is one distinguished maximum which is significantly larger from the rest of the values, which makes the lumping approximation more precise. However, the unscaled variant 1. and variant 2. are not appropriate for the lumping, since the large values could easily cause overflow. Therefore, when the large value is encountered the variant 3. is chosen. The disadvantage of variant 3. is that it causes more objective function and inequality constraints evaluations. For this reason, if the maximum value of the scaled constraints (variant 2) is below certain threshold we do not perform additional scaling by maximum as in variant 3. The overall algorithm can be summarized as follows:

\begin{lstlisting}
smoothing = 0.1
if (max(h(x)/f_max^2) > 5.0)
	alpha = smoothing * max(h(x)/f_max^2)
else
	alpha = smoothing
\end{lstlisting} 
The benchmarks such as case118, or case2737sop do converge, however PEGASE benchmarks do not. This is true also for the algorithm where no lumping is performed, therefore the convergence problem lies either in the Hessian approximation or the reduced space method itself (maybe insufficient precision of the PF solution).

\subsection*{Full-Space Quasi-Newton Convergence Study}
\begin{figure}[!ht]
	\centering
	\begin{subfigure}[b]{\columnwidth}
\definecolor{mycolor0}{HTML}{2B83BA}
\definecolor{mycolor1}{HTML}{D7191C}
\definecolor{mycolor2}{HTML}{FDAE61}
\definecolor{mycolor3}{HTML}{ABDDA4}

\begin{tikzpicture}

\pgfplotstableread{
iter   objective infpr    infdu
0   3.8755581000E+05    5.3400E+00  1.0000E+02
1   3.8735656000E+05    4.6300E+00  8.8200E+01
2   3.8713144000E+05    3.8100E+00  7.4200E+01
3   3.8709240000E+05    3.6700E+00  7.2400E+01
4   3.8705443000E+05    3.5300E+00  4.9300E+02
5   3.8671046000E+05    2.2500E+00  2.4400E+02
6   3.8657651000E+05    1.7500E+00  3.2800E+02
7   3.8638459000E+05    1.0300E+00  2.0700E+02
8   3.8631981000E+05    7.9300E-01  1.6700E+02
9   3.8623216000E+05    4.6700E-01  9.8000E+01
10  3.8619011000E+05    3.1100E-01  6.6500E+01
11  3.8617296000E+05    2.4700E-01  5.4600E+01
12  3.8614711000E+05    1.5100E-01  3.3900E+01
13  3.8613207000E+05    9.4700E-02  2.1500E+01
14  3.8612361000E+05    6.3200E-02  5.3700E+01
15  3.8612186000E+05    5.6700E-02  8.6300E+01
16  3.8612116000E+05    5.4100E-02  8.2600E+01
17  3.8612099000E+05    5.3500E-02  8.1700E+01
18  3.8611913000E+05    4.6600E-02  7.2300E+01
19  3.8611155000E+05    1.8400E-02  2.7100E+01
20  3.8611149000E+05    1.8100E-02  2.0900E+01
21  3.8611129000E+05    1.7400E-02  5.2900E+01
22  3.8610894000E+05    8.6800E-03  1.2000E+01
23  3.8610829000E+05    6.2300E-03  8.6500E+00
24  3.8610827000E+05    6.1800E-03  3.7200E+01
25  3.8610794000E+05    4.9500E-03  3.5500E+01
26  3.8610740000E+05    3.6100E-03  1.4700E+01
27  3.8610739000E+05    3.6500E-03  1.7100E+01
28  3.8610721000E+05    3.4100E-03  2.2300E+01
29  3.8610717000E+05    3.2400E-03  3.4100E+01
30  3.8610692000E+05    3.3500E-03  2.4700E+01
31  3.8610679000E+05    2.4400E-03  1.9400E+01
32  3.8610672000E+05    1.5000E-03  8.9200E+00
33  3.8610669000E+05    1.2000E-02  2.1300E+01
34  3.8610669000E+05    1.1300E-02  1.4400E+01
35  3.8610669000E+05    1.1200E-02  1.8300E+01
36  3.8610668000E+05    1.2900E-02  5.5000E+01
37  3.8610665000E+05    7.8700E-03  3.3400E+01
38  3.8610663000E+05    5.3700E-03  2.7800E+01
39  3.8610663000E+05    3.2000E-02  2.4600E+01
40  3.8610663000E+05    2.3200E-02  2.2500E+01
41  3.8610662000E+05    3.7700E-02  1.9800E+01
42  3.8610662000E+05    4.8800E-02  1.6300E+01
43  3.8610661000E+05    8.7800E-02  2.0000E+00
44  3.8610661000E+05    4.3900E-02  2.9200E-03
45  3.8610661000E+05    1.4600E-02  1.4800E+00
46  3.8610661000E+05    1.1000E-03  3.1700E-03
}\datatable;


\begin{axis}[
width=0.97\linewidth, height=3.5cm,
xmin=0,
xmax=50,
axis lines*=left,
yticklabel style={
	/pgf/number format/fixed,
	/pgf/number format/fixed zerofill,
	/pgf/number format/precision=2
},
xtick={0,10,20,30,40,50,60,70,80,90},
yticklabel style={font=\scriptsize},
xticklabel style={font=\scriptsize},
ylabel style={font=\scriptsize},
xlabel style={font=\footnotesize},
ylabel={Objective Value (\$/h)},
y label style={at={(axis description cs:-0.08,.5)},anchor=south},
xlabel={IPM Iteration number},
legend=none    
]
\addplot[color=mycolor0, thick] table[x=iter, y=objective] {\datatable}; \label{plot:objE13k}
\end{axis}

\begin{axis}[name=symb,
width=0.97\linewidth, height=3.5cm,
ymode=log, 
log origin=infty,
axis lines*=right, 
axis x line=none, 
ymajorgrids, yminorgrids,
xmin=0,
xmax=50,
ytick={1e-5, 1e-3, 1e-2, 1e-1, 1e0, 1e1, 1e2, 1e3},
xticklabel style={rotate=0, xshift=-0.0cm, anchor=north, font=\scriptsize},
yticklabel style={font=\scriptsize},
ylabel style={font=\scriptsize},
xlabel style={font=\scriptsize},
ylabel={Feasibility},
y label style={at={(axis description cs:1.0,.5)},anchor=south},
legend style={at={(0.9,1.06)},legend cell align=left,align=right,draw=none,fill=none,font=\tiny,legend columns=3}
]
\addlegendimage{/pgfplots/refstyle=plot:objE13k}\addlegendentry{Objective}

\addplot[color=mycolor1, thick] table[x=iter, y=infpr] {\datatable};
\addlegendentry{Primal Feasibility}


\addplot[color=mycolor2, thick] table[x=iter, y=infdu] {\datatable}; 
\addlegendentry{Dual Feasibility}
  
\end{axis}

\end{tikzpicture}
	\vspace{-0.5cm}
	\caption{Exact Newton. \label{fig:convergence13kExact}}
	\end{subfigure}
	\begin{subfigure}[b]{\columnwidth}
	\input{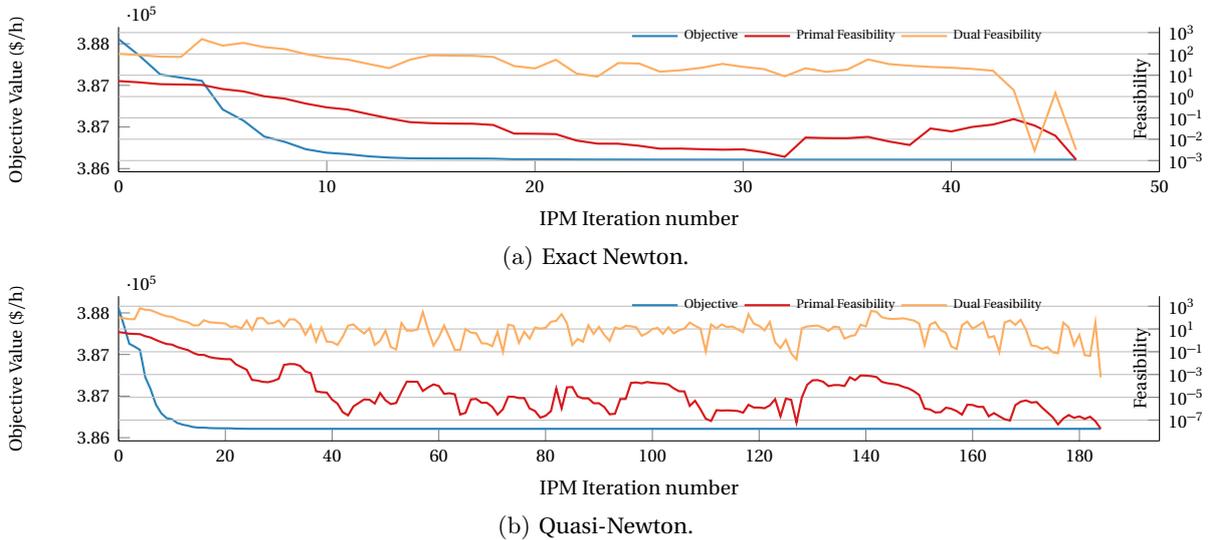}
	\vspace{-0.5cm}
	\caption{Quasi-Newton. \label{fig:convergence13kInexact}}
	\end{subfigure}
	\caption{Convergence trajectory for PEGASE13659-10 power grid benchmark, considering the full space IPM. Tolerance fixed at $tol = 10^{-2}$. \label{fig:convergence13k}}
\end{figure}

It is well known that quasi-Newton's methods have worse convergence that their exact counterparts. The convergence of the quasi-Newton's full space approach is studied in this section. Figures \ref{fig:convergence13k}
and \ref{fig:convergence1k} 
illustrate the convergence of both exact and quasi-Newton full-space approaches for two benchmarks. In both cases, the quasi-Newton requires significantly more iterations to reach the required tolerance $\epsilon_{tol} = 10^{-2}$. However, it is important to note that the quasi-Newton method reaches the optimal value and constraint (primal) feasibility within a similar number of iterations as required by the exact Newton method. The remaining iterations are performed in order to satisfy the tolerance for the dual feasibility (i.e., optimality). We can argue that since we have observed that the objective function is optimal and constraint feasibility is satisfied, the IPM iterations can be stopped early by relaxing the dual infeasibility tolerance.  The convergence problems might occur due to the fact that the BFGS cannot approximate the Hessian at late iterates (close to the optimal point) due to the well-known ill-conditioning of the Hessian in the IPM methods as the solution is approached. Exactly for this reason, since we are close to the solution, we can terminate early, not waiting for the dual feasibility reaching the desired tolerance.
It is also noted in the IPOPT manual\footnote{https://www.coin-or.org/Ipopt/documentation/node42.html} that quasi-Newton methods have trouble bringing down the dual infeasibility. Therefore, special acceptable convergence tolerance is provided. If the usual tolerances are satisfied at an iteration, the algorithm immediately terminates with a success message. On the other hand, if the algorithm encounters many iterations in a row that are considered acceptable, it will terminate before the desired convergence tolerance is met. Such acceptance stopping criterion may be based on the objective function change.



\begin{figure}[!ht]
	\centering
	\begin{subfigure}[b]{\columnwidth}
\definecolor{mycolor0}{HTML}{2B83BA}
\definecolor{mycolor1}{HTML}{D7191C}
\definecolor{mycolor2}{HTML}{FDAE61}
\definecolor{mycolor3}{HTML}{ABDDA4}

\begin{tikzpicture}

\pgfplotstableread{
iter   objective infpr    infdu
0   7.4014002E+04   2.670E+00   1.000E+02
1   7.4019744E+04   2.310E+00   8.770E+01
2   7.4027313E+04   1.900E+00   7.280E+01
3   7.4031228E+04   1.710E+00   6.680E+01
4   7.4039160E+04   1.350E+00   5.420E+01
5   7.4044578E+04   1.110E+00   4.430E+01
6   7.4051439E+04   8.030E-01   3.210E+01
7   7.4054128E+04   6.840E-01   2.510E+02
8   7.4059307E+04   4.530E-01   1.970E+02
9   7.4060179E+04   4.140E-01   1.940E+02
10  7.4060226E+04   4.120E-01   2.210E+02
11  7.4062611E+04   3.030E-01   1.650E+02
12  7.4062732E+04   2.970E-01   1.660E+02
13  7.4065734E+04   1.620E-01   8.960E+01
14  7.4067516E+04   8.290E-02   4.560E+01
15  7.4067783E+04   7.100E-02   3.990E+01
16  7.4068231E+04   5.100E-02   2.880E+01
17  7.4069029E+04   1.580E-02   2.200E+01
18  7.4069387E+04   3.500E-03   8.900E-03
}\datatable;


\begin{axis}[
width=0.97\linewidth, height=3.5cm,
xmin=0,
xmax=20,
axis lines*=left,
yticklabel style={
	/pgf/number format/fixed,
	/pgf/number format/fixed zerofill,
	/pgf/number format/precision=2
},
xtick={0,5,10,15,20,30,40,50,60,70,80,90},
yticklabel style={font=\scriptsize},
xticklabel style={font=\scriptsize},
ylabel style={font=\scriptsize},
xlabel style={font=\footnotesize},
ylabel={Objective Value (\$/h)},
y label style={at={(axis description cs:-0.08,.5)},anchor=south},
xlabel={IPM Iteration number},
legend=none    
]
\addplot[color=mycolor0, thick] table[x=iter, y=objective] {\datatable}; \label{plot:objE1k}
\end{axis}

\begin{axis}[name=symb,
width=0.97\linewidth, height=3.5cm,
ymode=log, 
log origin=infty,
axis lines*=right, 
axis x line=none, 
ymajorgrids, yminorgrids,
xmin=0,
xmax=20,
ytick={1e-5, 1e-3, 1e-2, 1e-1, 1e0, 1e1, 1e2, 1e3},
xticklabel style={rotate=0, xshift=-0.0cm, anchor=north, font=\scriptsize},
yticklabel style={font=\scriptsize},
ylabel style={font=\scriptsize},
xlabel style={font=\scriptsize},
ylabel={Feasibility},
y label style={at={(axis description cs:1.0,.5)},anchor=south},
legend style={at={(1.0,1.25)},legend cell align=left,align=right,draw=none,fill=none,font=\tiny,legend columns=3}
]
\addlegendimage{/pgfplots/refstyle=plot:objE1k}\addlegendentry{Objective}

\addplot[color=mycolor1, thick] table[x=iter, y=infpr] {\datatable};
\addlegendentry{Primal Feasibility}


\addplot[color=mycolor2, thick] table[x=iter, y=infdu] {\datatable}; 
\addlegendentry{Dual Feasibility}
  
\end{axis}

\end{tikzpicture}
	\vspace{-0.5cm}
	\caption{Exact Newton. \label{fig:convergence1kExact}}
	\end{subfigure}
	\begin{subfigure}[b]{\columnwidth}
	\input{figures/fig_convergence1354_inexact.tex}
	\vspace{-0.5cm}
	\caption{Quasi-Newton. \label{fig:convergence1kInexact}}
	\end{subfigure}
	\caption{Convergence trajectory for PEGASE1354-5 power grid benchmark, considering the full space IPM. Tolerance fixed at $tol = 10^{-2}$. \label{fig:convergence1k}}
\end{figure}




\begin{figure}[!ht]
	\centering
	\begin{subfigure}[b]{\columnwidth}
\definecolor{mycolor0}{HTML}{2B83BA}
\definecolor{mycolor1}{HTML}{D7191C}
\definecolor{mycolor2}{HTML}{FDAE61}
\definecolor{mycolor3}{HTML}{ABDDA4}

\begin{tikzpicture}

\pgfplotstableread{
iter   objective infpr    infdu
0	1.33E+05	1.43E+01	4.53E+00
1	1.33E+05	1.36E+01	4.26E+00
2	1.32E+05	1.14E+01	7.22E+00
3	1.31E+05	7.73E+00	2.20E+01
4	1.31E+05	3.82E+00	4.32E+01
5	1.30E+05	2.37E+00	3.04E+01
6	1.31E+05	2.00E-01	4.88E+01
7	1.31E+05	1.04E-02	2.06E+00
8	1.30E+05	1.20E-03	3.00E+00
9	1.30E+05	7.29E-04	4.64E+00
10	1.30E+05	5.60E-04	1.68E+00
11	1.30E+05	1.33E-04	3.33E-02
12	1.30E+05	1.13E-05	1.19E-03
13	1.30E+05	1.38E-06	9.18E-05
}\datatable;


\begin{axis}[
width=0.97\linewidth, height=3.5cm,
xmin=0,
xmax=15,
axis lines*=left,
yticklabel style={
	/pgf/number format/fixed,
	/pgf/number format/fixed zerofill,
	/pgf/number format/precision=2
},
yticklabel style={font=\scriptsize},
xticklabel style={font=\scriptsize},
ylabel style={font=\scriptsize},
xlabel style={font=\footnotesize},
ylabel={Objective Value (\$/h)},
y label style={at={(axis description cs:-0.08,.5)},anchor=south},
xlabel={IPM Iteration number},
legend=none    
]
\addplot[color=mycolor0, thick] table[x=iter, y=objective] {\datatable}; \label{plot:objE1k}
\end{axis}

\begin{axis}[name=symb,
width=0.97\linewidth, height=3.5cm,
ymode=log, 
log origin=infty,
axis lines*=right, 
axis x line=none, 
ymajorgrids,
xmin=0,
xmax=15,
ytick={1e-6,1e-5,1e-4, 1e-3, 1e-2, 1e-1, 1e0, 1e1, 1e2, 1e3},
xticklabel style={rotate=0, xshift=-0.0cm, anchor=north, font=\scriptsize},
yticklabel style={font=\scriptsize},
ylabel style={font=\scriptsize},
xlabel style={font=\scriptsize},
ylabel={Feasibility},
y label style={at={(axis description cs:1.0,.5)},anchor=south},
legend style={at={(1.0,1.25)},legend cell align=left,align=right,draw=none,fill=none,font=\tiny,legend columns=3}
]
\addlegendimage{/pgfplots/refstyle=plot:objE1k}\addlegendentry{Objective}

\addplot[color=mycolor1, thick] table[x=iter, y=infpr] {\datatable};
\addlegendentry{Primal Feasibility}


\addplot[color=mycolor2, thick] table[x=iter, y=infdu] {\datatable}; 
\addlegendentry{Dual Feasibility}
  
\end{axis}

\end{tikzpicture}
	\vspace{-0.5cm}
	\caption{Exact Newton. \label{fig:convergence118Exact}}
	\end{subfigure}
	\begin{subfigure}[b]{\columnwidth}
\definecolor{mycolor0}{HTML}{2B83BA}
\definecolor{mycolor1}{HTML}{D7191C}
\definecolor{mycolor2}{HTML}{FDAE61}
\definecolor{mycolor3}{HTML}{ABDDA4}

\begin{tikzpicture}

\pgfplotstableread{
iter   objective infpr    infdu
0	1.33E+05	1.43E+01	4.53E+00
1	1.33E+05	1.37E+01	1.00E+01
2	1.32E+05	1.15E+01	5.29E+01
3	1.31E+05	7.83E+00	4.55E+01
4	1.31E+05	5.33E+00	7.07E+01
5	1.31E+05	3.98E+00	6.69E+01
6	1.31E+05	1.61E+00	4.94E+01
7	1.31E+05	1.45E+00	3.63E+01
8	1.31E+05	5.02E-03	1.12E+02
9	1.31E+05	4.09E-02	3.70E+01
10	1.31E+05	3.27E-02	5.47E+01
11	1.31E+05	1.03E-04	1.51E+01
12	1.31E+05	8.53E-04	4.68E+00
13	1.31E+05	2.25E-03	5.65E+01
14	1.31E+05	2.75E-03	4.41E+01
15	1.31E+05	1.05E-03	1.94E+02
16	1.31E+05	4.68E-03	1.10E+01
17	1.31E+05	2.82E-03	5.79E+00
18	1.31E+05	2.31E-03	3.41E+01
19	1.31E+05	1.80E-03	4.40E+01
20	1.31E+05	2.55E-04	7.37E+00
21	1.31E+05	1.10E-04	8.69E+00
22	1.31E+05	6.42E-05	1.04E+01
23	1.31E+05	4.50E-05	6.13E+00
24	1.31E+05	1.55E-04	1.95E+01
25	1.31E+05	5.77E-04	1.24E+01
26	1.31E+05	4.00E-04	2.87E+01
27	1.31E+05	3.78E-03	3.51E+01
28	1.31E+05	1.15E-03	4.42E+01
29	1.31E+05	3.58E-03	1.04E+01
30	1.31E+05	3.48E-04	7.45E+00
31	1.31E+05	1.85E-04	4.27E+01
32	1.31E+05	3.58E-04	3.68E+00
33	1.31E+05	4.43E-05	6.68E+00
34	1.31E+05	1.28E-04	4.74E+00
35	1.31E+05	1.37E-04	6.14E+00
36	1.31E+05	1.90E-04	6.38E+00
37	1.31E+05	3.69E-04	1.44E+01
38	1.30E+05	3.40E-04	1.53E+01
39	1.31E+05	2.77E-04	4.53E+01
40	1.31E+05	5.57E-04	1.33E+00
41	1.31E+05	4.94E-04	1.07E+01
42	1.30E+05	4.29E-04	1.13E+01
43	1.31E+05	1.43E-04	2.23E+01
44	1.31E+05	1.39E-04	2.88E+01
45	1.30E+05	1.12E-04	1.68E+01
46	1.30E+05	2.64E-05	2.34E+00
47	1.30E+05	1.95E-05	9.36E+00
48	1.30E+05	1.95E-05	6.92E+00
49	1.30E+05	7.11E-05	5.42E+00
50	1.30E+05	7.38E-05	1.56E+01
51	1.30E+05	9.01E-05	1.68E+01
52	1.30E+05	1.34E-04	6.62E+00
53	1.30E+05	2.22E-04	8.45E+00
54	1.30E+05	3.30E-04	1.08E+01
55	1.30E+05	3.29E-04	1.12E+01
56	1.30E+05	3.30E-04	1.44E+01
57	1.30E+05	3.11E-04	5.27E+00
58	1.31E+05	4.58E-05	1.36E+01
59	1.30E+05	1.09E-03	2.31E+00
60	1.30E+05	9.82E-04	8.63E+00
61	1.31E+05	2.71E-04	2.01E+01
62	1.30E+05	2.91E-04	9.05E+00
63	1.31E+05	1.03E-05	2.58E+01
64	1.30E+05	3.71E-04	2.48E-01
65	1.30E+05	6.67E-05	9.93E-01
66	1.30E+05	6.33E-05	7.82E+00
67	1.30E+05	4.93E-05	2.21E+00
68	1.31E+05	1.09E-04	2.53E+01
69	1.31E+05	5.94E-04	3.36E+01
70	1.31E+05	5.88E-04	2.99E+01
71	1.30E+05	3.68E-03	1.52E+01
72	1.30E+05	3.61E-03	1.15E+01
73	1.30E+05	5.10E-04	4.91E+00
74	1.31E+05	4.56E-05	1.02E+02
75	1.31E+05	8.89E-05	3.58E+01
76	1.31E+05	2.54E-04	3.62E+01
77	1.31E+05	4.27E-04	3.69E+01
78	1.31E+05	4.29E-04	3.00E+01
79	1.31E+05	4.39E-04	2.78E+01
80	1.31E+05	1.91E-03	1.62E+01
81	1.31E+05	3.50E-04	5.34E+01
82	1.31E+05	1.50E-03	4.58E+01
83	1.31E+05	1.27E-03	5.23E+01
84	1.31E+05	1.34E-05	6.26E+01
85	1.31E+05	4.10E-03	1.99E+01
86	1.30E+05	3.50E-03	1.26E+01
87	1.30E+05	1.39E-04	1.09E+00
88	1.30E+05	1.01E-05	2.79E+00
89	1.30E+05	1.18E-05	5.97E-01
90	1.30E+05	6.61E-06	4.51E+00
91	1.30E+05	3.23E-05	4.50E+00
92	1.30E+05	4.90E-05	5.58E+00
93	1.30E+05	4.20E-05	3.56E+00
94	1.30E+05	1.93E-07	3.50E+00
95	1.30E+05	2.26E-06	3.86E+00
96	1.30E+05	1.39E-05	8.66E-01
97	1.30E+05	6.57E-06	1.06E+00
98	1.30E+05	5.67E-06	6.92E-01
99	1.30E+05	5.31E-06	9.39E-01
100	1.30E+05	9.12E-06	3.31E+00
101	1.30E+05	8.58E-06	8.64E+00
102	1.30E+05	1.60E-05	1.03E+01
103	1.30E+05	3.71E-05	3.00E+00
104	1.30E+05	1.85E-05	8.14E-01
105	1.30E+05	3.50E-06	6.53E-01
106	1.30E+05	1.08E-06	2.50E-01
107	1.30E+05	1.26E-05	6.66E-01
108	1.30E+05	1.23E-05	1.49E+00
109	1.30E+05	1.27E-05	1.86E+00
110	1.30E+05	1.33E-05	3.12E+00
111	1.30E+05	7.69E-06	9.83E+00
112	1.30E+05	7.36E-06	9.42E+00
113	1.30E+05	7.36E-06	9.42E+00
114	1.30E+05	2.36E-05	1.01E+00
115	1.30E+05	7.33E-07	8.02E-02
116	1.30E+05	1.31E-07	2.37E-01
117	1.30E+05	2.49E-06	1.20E+00
118	1.30E+05	2.49E-06	1.69E+00
119	1.30E+05	3.40E-05	1.48E+00
120	1.30E+05	3.29E-05	1.35E+00
121	1.30E+05	3.08E-05	1.22E+00
122	1.30E+05	3.00E-05	2.38E+00
123	1.30E+05	1.02E-05	1.48E+00
124	1.30E+05	1.05E-05	7.84E-01
125	1.30E+05	3.10E-06	4.24E-01
126	1.30E+05	1.52E-06	2.24E+00
127	1.30E+05	1.35E-06	2.08E-01
128	1.30E+05	2.97E-07	1.41E-01
129	1.30E+05	4.70E-07	7.37E-01
130	1.30E+05	3.01E-06	6.96E-01
131	1.30E+05	1.75E-05	1.85E+00
132	1.30E+05	2.32E-05	2.63E+00
133	1.30E+05	3.62E-05	3.26E+00
134	1.30E+05	3.84E-05	3.79E+00
135	1.30E+05	3.25E-05	3.63E+00
136	1.30E+05	3.22E-05	2.87E+00
137	1.30E+05	3.78E-05	3.70E+00
138	1.30E+05	6.13E-06	7.86E+00
139	1.30E+05	1.05E-05	2.07E+00
140	1.30E+05	8.89E-06	8.30E-01
141	1.30E+05	1.87E-08	5.78E+00
142	1.30E+05	9.63E-06	1.88E+00
143	1.30E+05	3.99E-05	1.60E+00
144	1.30E+05	1.68E-05	2.87E+00
145	1.30E+05	8.58E-06	2.98E+00
146	1.30E+05	4.06E-06	1.66E+00
147	1.30E+05	9.90E-06	4.82E+00
148	1.30E+05	5.75E-06	9.28E-02
149	1.30E+05	1.30E-08	4.23E-02
150	1.30E+05	2.97E-07	9.46E-02
}\datatable;


\begin{axis}[
width=0.97\linewidth, height=3.5cm,
xmin=0,
xmax=155,
axis lines*=left,
yticklabel style={
	/pgf/number format/fixed,
	/pgf/number format/fixed zerofill,
	/pgf/number format/precision=2
},
yticklabel style={font=\scriptsize},
xticklabel style={font=\scriptsize},
ylabel style={font=\scriptsize},
xlabel style={font=\footnotesize},
ylabel={Objective Value (\$/h)},
y label style={at={(axis description cs:-0.08,.5)},anchor=south},
xlabel={IPM Iteration number},
legend=none    
]
\addplot[color=mycolor0, thick] table[x=iter, y=objective] {\datatable}; \label{plot:objE1k}
\end{axis}

\begin{axis}[name=symb,
width=0.97\linewidth, height=3.5cm,
ymode=log, 
log origin=infty,
axis lines*=right, 
axis x line=none, 
ymajorgrids,
xmin=0,
xmax=155,
ytick={1e-6,1e-5,1e-4, 1e-3, 1e-2, 1e-1, 1e0, 1e1, 1e2, 1e3},
xticklabel style={rotate=0, xshift=-0.0cm, anchor=north, font=\scriptsize},
yticklabel style={font=\scriptsize},
ylabel style={font=\scriptsize},
xlabel style={font=\scriptsize},
ylabel={Feasibility},
y label style={at={(axis description cs:1.0,.5)},anchor=south},
legend style={at={(1.0,1.25)},legend cell align=left,align=right,draw=none,fill=none,font=\tiny,legend columns=3}
]
\addlegendimage{/pgfplots/refstyle=plot:objE1k}\addlegendentry{Objective}

\addplot[color=mycolor1, thick] table[x=iter, y=infpr] {\datatable};
\addlegendentry{Primal Feasibility}


\addplot[color=mycolor2, thick] table[x=iter, y=infdu] {\datatable}; 
\addlegendentry{Dual Feasibility}
  
\end{axis}

\end{tikzpicture}
	\vspace{-0.5cm}
	\caption{Full-Space Quasi-Newton. \label{fig:convergence118Inexact}}
	\end{subfigure}
	\begin{subfigure}[b]{\columnwidth}
\definecolor{mycolor0}{HTML}{2B83BA}
\definecolor{mycolor1}{HTML}{D7191C}
\definecolor{mycolor2}{HTML}{FDAE61}
\definecolor{mycolor3}{HTML}{ABDDA4}

\begin{tikzpicture}

\pgfplotstableread{
iter   objective infpr    infdu
0	1.31E+05	1.43E+01	9.98E+00
1	1.32E+05	1.36E+01	3.25E+02
2	1.31E+05	1.13E+01	4.33E+02
3	1.31E+05	8.73E+00	2.30E+02
4	1.31E+05	7.62E+00	1.95E+02
5	1.31E+05	5.55E+00	1.06E+02
6	1.32E+05	3.65E+00	3.95E+02
7	1.32E+05	1.86E+01	4.88E+02
8	1.33E+05	9.20E+00	9.63E+02
9	1.33E+05	5.01E+00	8.94E+02
10	1.33E+05	8.12E+00	1.12E+03
11	1.33E+05	0.00E+00	6.67E+02
12	1.32E+05	0.00E+00	1.01E+02
13	1.32E+05	0.00E+00	4.13E+01
14	1.32E+05	0.00E+00	1.15E+02
15	1.32E+05	0.00E+00	3.61E+02
16	1.32E+05	1.83E+00	6.69E+02
17	1.31E+05	0.00E+00	2.19E+02
18	1.32E+05	1.27E+01	9.19E+02
19	1.31E+05	0.00E+00	2.18E+02
20	1.31E+05	1.92E+00	4.60E+02
21	1.31E+05	0.00E+00	8.36E+01
22	1.31E+05	0.00E+00	2.24E+02
23	1.31E+05	3.09E-03	1.36E+02
24	1.31E+05	0.00E+00	1.81E+02
25	1.31E+05	0.00E+00	1.56E+02
26	1.31E+05	0.00E+00	3.29E+02
27	1.31E+05	0.00E+00	3.54E+01
28	1.31E+05	0.00E+00	1.78E+02
29	1.31E+05	0.00E+00	8.05E+01
30	1.31E+05	0.00E+00	6.69E+01
31	1.31E+05	0.00E+00	6.74E+01
32	1.31E+05	0.00E+00	2.08E+01
33	1.31E+05	0.00E+00	4.18E+01
34	1.31E+05	0.00E+00	4.86E+01
35	1.31E+05	4.82E-01	7.46E+01
36	1.31E+05	1.18E-02	2.32E+02
37	1.31E+05	0.00E+00	1.80E+01
38	1.31E+05	0.00E+00	8.64E+00
39	1.31E+05	0.00E+00	3.46E+01
40	1.31E+05	0.00E+00	2.03E+01
41	1.31E+05	0.00E+00	4.43E+01
42	1.31E+05	0.00E+00	2.37E+01
43	1.31E+05	0.00E+00	1.11E+02
44	1.31E+05	0.00E+00	1.78E+02
45	1.31E+05	0.00E+00	1.21E+02
46	1.31E+05	0.00E+00	6.42E+01
47	1.31E+05	0.00E+00	7.77E+00
48	1.31E+05	0.00E+00	3.04E+01
49	1.31E+05	0.00E+00	1.43E+01
50	1.31E+05	0.00E+00	3.93E+01
51	1.31E+05	0.00E+00	3.68E+01
52	1.31E+05	0.00E+00	2.18E+01
53	1.31E+05	0.00E+00	1.24E+01
54	1.30E+05	0.00E+00	1.06E+01
55	1.30E+05	0.00E+00	2.14E+01
56	1.30E+05	0.00E+00	3.36E+01
57	1.31E+05	0.00E+00	1.23E+02
58	1.31E+05	0.00E+00	2.71E+01
59	1.31E+05	0.00E+00	3.64E+01
60	1.30E+05	0.00E+00	6.38E+01
61	1.30E+05	0.00E+00	2.74E+01
62	1.30E+05	0.00E+00	2.54E+01
63	1.30E+05	0.00E+00	1.01E+01
64	1.30E+05	0.00E+00	4.02E+01
65	1.30E+05	3.47E-03	2.05E+01
66	1.30E+05	0.00E+00	4.09E+01
67	1.30E+05	0.00E+00	4.66E+01
68	1.30E+05	0.00E+00	7.59E+01
69	1.30E+05	0.00E+00	6.70E+00
70	1.30E+05	0.00E+00	3.38E+00
71	1.30E+05	0.00E+00	5.62E+00
72	1.30E+05	0.00E+00	1.21E+01
73	1.30E+05	0.00E+00	1.11E+01
74	1.30E+05	0.00E+00	1.25E+01
75	1.30E+05	0.00E+00	9.81E+00
76	1.30E+05	0.00E+00	2.81E+01
77	1.31E+05	0.00E+00	1.32E+02
78	1.30E+05	0.00E+00	1.13E+01
79	1.30E+05	0.00E+00	2.43E+01
80	1.30E+05	0.00E+00	4.08E+00
81	1.30E+05	0.00E+00	3.45E+00
82	1.30E+05	0.00E+00	2.61E+00
83	1.30E+05	0.00E+00	3.32E+00
84	1.30E+05	0.00E+00	6.74E+00
85	1.30E+05	0.00E+00	4.28E+00
86	1.30E+05	0.00E+00	1.48E+01
87	1.30E+05	0.00E+00	1.62E+01
88	1.30E+05	0.00E+00	6.12E+00
89	1.30E+05	0.00E+00	4.37E+00
90	1.30E+05	0.00E+00	3.37E+00
91	1.30E+05	0.00E+00	4.71E+00
92	1.30E+05	0.00E+00	1.22E+01
93	1.30E+05	0.00E+00	9.34E+00
94	1.30E+05	0.00E+00	3.23E+01
95	1.30E+05	0.00E+00	7.05E+00
96	1.30E+05	0.00E+00	2.07E+01
97	1.30E+05	0.00E+00	5.55E-01
98	1.30E+05	0.00E+00	5.12E-01
99	1.30E+05	0.00E+00	4.61E+00
100	1.30E+05	0.00E+00	3.47E+00
101	1.30E+05	5.93E-07	3.14E+00
102	1.30E+05	4.73E-07	2.51E+00
103	1.30E+05	0.00E+00	1.27E+01
104	1.30E+05	0.00E+00	8.21E+00
105	1.30E+05	0.00E+00	5.84E+01
106	1.30E+05	0.00E+00	5.84E+00
107	1.30E+05	0.00E+00	1.76E+01
108	1.30E+05	0.00E+00	2.25E+01
109	1.30E+05	0.00E+00	5.35E+00
110	1.30E+05	0.00E+00	7.28E+00
111	1.30E+05	0.00E+00	6.68E+00
112	1.30E+05	0.00E+00	1.61E+00
113	1.30E+05	0.00E+00	2.85E+00
114	1.30E+05	0.00E+00	6.40E-01
115	1.30E+05	0.00E+00	2.31E+00
116	1.30E+05	0.00E+00	1.67E+00
117	1.30E+05	0.00E+00	5.95E-01
118	1.30E+05	0.00E+00	2.76E-01
119	1.30E+05	0.00E+00	3.69E-01
120	1.30E+05	0.00E+00	6.74E-01
121	1.30E+05	0.00E+00	1.65E+00
122	1.30E+05	0.00E+00	2.08E+00
123	1.30E+05	0.00E+00	1.92E+01
124	1.30E+05	0.00E+00	2.56E+01
125	1.30E+05	0.00E+00	7.80E+00
126	1.30E+05	0.00E+00	9.61E+00
127	1.30E+05	0.00E+00	1.83E+00
128	1.30E+05	0.00E+00	1.03E+01
129	1.30E+05	0.00E+00	8.40E-01
130	1.30E+05	0.00E+00	3.60E+00
131	1.30E+05	0.00E+00	5.14E-01
132	1.30E+05	0.00E+00	1.06E+00
133	1.30E+05	0.00E+00	4.79E-01
134	1.30E+05	0.00E+00	2.96E-01
135	1.30E+05	0.00E+00	7.04E-01
136	1.30E+05	0.00E+00	7.92E-01
137	1.30E+05	0.00E+00	1.62E+00
138	1.30E+05	0.00E+00	1.54E-01
139	1.30E+05	0.00E+00	3.10E-01
140	1.30E+05	0.00E+00	1.16E-01
141	1.30E+05	0.00E+00	2.57E-01
142	1.30E+05	0.00E+00	2.14E-01
143	1.30E+05	4.85E-07	6.39E-01
144	1.30E+05	1.62E-06	3.90E-01
145	1.30E+05	0.00E+00	9.21E+00
146	1.30E+05	7.41E-05	1.39E+00
147	1.30E+05	8.88E-09	5.36E+00
148	1.30E+05	8.79E-05	1.52E+00
149	1.30E+05	0.00E+00	1.97E+01
150	1.30E+05	2.93E-07	8.99E+00
151	1.30E+05	2.90E-05	1.04E+00
152	1.30E+05	1.45E-06	2.65E-01
153	1.30E+05	0.00E+00	9.23E-02
154	1.30E+05	7.64E-07	1.34E-01
155	1.30E+05	5.57E-08	7.75E-02
156	1.30E+05	5.28E-09	4.94E-02
157	1.30E+05	0.00E+00	9.81E-02
}\datatable;


\begin{axis}[
width=0.97\linewidth, height=3.5cm,
xmin=0,
xmax=160,
axis lines*=left,
yticklabel style={
	/pgf/number format/fixed,
	/pgf/number format/fixed zerofill,
	/pgf/number format/precision=2
},
yticklabel style={font=\scriptsize},
xticklabel style={font=\scriptsize},
ylabel style={font=\scriptsize},
xlabel style={font=\footnotesize},
ylabel={Objective Value (\$/h)},
y label style={at={(axis description cs:-0.08,.5)},anchor=south},
xlabel={IPM Iteration number},
legend=none    
]
\addplot[color=mycolor0, thick] table[x=iter, y=objective] {\datatable}; \label{plot:objE1k}
\end{axis}

\begin{axis}[name=symb,
width=0.97\linewidth, height=3.5cm,
ymode=log, 
log origin=infty,
axis lines*=right, 
axis x line=none, 
ymajorgrids,
xmin=0,
xmax=160,
ytick={1e-6,1e-5, 1e-4,1e-3, 1e-2, 1e-1, 1e0, 1e1, 1e2, 1e3},
xticklabel style={rotate=0, xshift=-0.0cm, anchor=north, font=\scriptsize},
yticklabel style={font=\scriptsize},
ylabel style={font=\scriptsize},
xlabel style={font=\scriptsize},
ylabel={Feasibility},
y label style={at={(axis description cs:1.0,.5)},anchor=south},
legend style={at={(1.0,1.25)},legend cell align=left,align=right,draw=none,fill=none,font=\tiny,legend columns=3}
]
\addlegendimage{/pgfplots/refstyle=plot:objE1k}\addlegendentry{Objective}

\addplot[color=mycolor1, thick] table[x=iter, y=infpr] {\datatable};
\addlegendentry{Primal Feasibility}


\addplot[color=mycolor2, thick] table[x=iter, y=infdu] {\datatable}; 
\addlegendentry{Dual Feasibility}
  
\end{axis}

\end{tikzpicture}
	\vspace{-0.5cm}
	\caption{Reduced-Space Quasi-Newton. \label{fig:convergence118Adjoint}}
	\end{subfigure}
	\caption{Convergence trajectory for case118 power grid benchmark. \label{fig:convergence118}}
\end{figure}
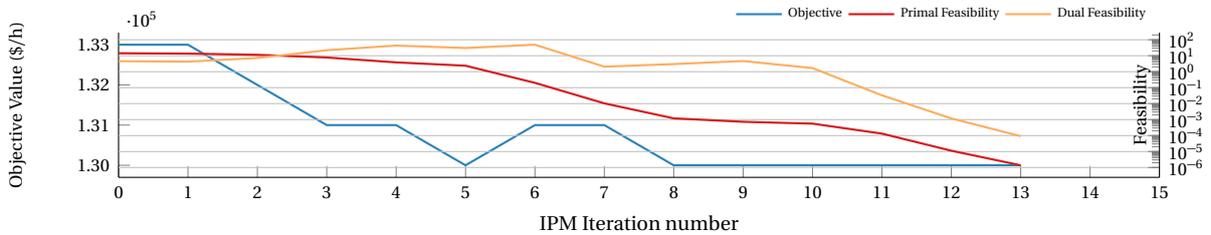
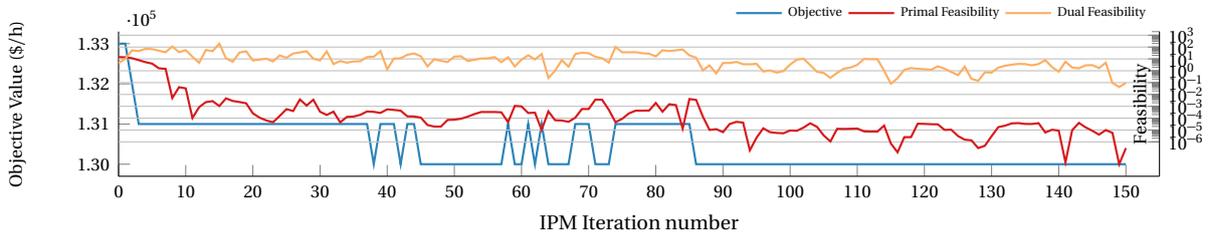
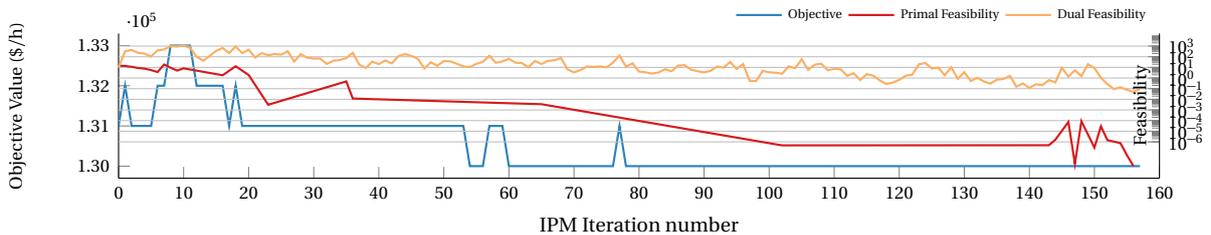

\section{Conclusions}
The optimization problems in the PDE-constrained problems (e.g., maximization of the oil production \cite{DrososConstraints}) are often solved by reduced space optimization methods \cite{Ghattas-01, Ghattas-02}. The discretized state variables and the PDEs representing the nonlinear equality constraints are removed from the optimization problem and are treated explicitly during the evaluation of the objective function value and its gradient. Given the control variables, the system of the removed equality constraints is used to solve for the state variables. Thus, the equality constraints are implicitly satisfied. The efficient evaluation of the gradient information is achieved using the adjoint method. The computational cost of evaluating the inequality constraints gradients can be reduced using the constraint lumping techniques \cite{DrososConstraints}. The second order derivatives are usually not evaluated exactly due to the excessive computational cost, only approximations such as BFGS are used. 

This approach was applied to the nonconvex OPF problems and the the results demonstrated that the computational cost of evaluating the gradient information is excessive and the constraint lumping introduces non-smooth functions, which leads to convergence difficulties of the IP method used for solution of the OPF problems. Additional investigation needs to be performed, closely analyzing the Hessian approximation and smoothing technique used during the constraints evaluation.



\bibliographystyle{plain}
\bibliography{biblio}

\end{document}